\documentclass[leqno,11pt]{amsart}
\usepackage{amssymb}
\usepackage{amsmath}
\usepackage{latexsym}%
\usepackage[normalem]{ulem}
\usepackage{color}
\usepackage[myheadings]{fullpage}
\usepackage[all]{xy}
\usepackage[linktoc=all, pagebackref, hyperindex]{hyperref}%
\hypersetup{colorlinks,  citecolor=blue,  filecolor=blue,  linkcolor=red,  urlcolor=black}
\usepackage{marginnote} 

\newcommand{\s}{\mathbf{s}}

\newcommand{\m}{\mathfrak{m}}
\newcommand{\n}{\mathfrak{n}}

\def\NN{{\mathbb{N}}}
\def\QQ{{\mathbb{Q}}}

\def\KK{{\mathbb{K}}}

\def\NN{{\mathbb{N}}}
\def\QQ{{\mathbb{Q}}}

\def\KK{{\mathbb{K}}}
\def\DD{{\mathbb{D}}}

\DeclareMathOperator{\slop}{slop}
\DeclareMathOperator{\edim}{edim}

\DeclareMathOperator{\Hom}{Hom}

\DeclareMathOperator{\Ap}{Ap}
\DeclareMathOperator{\tr}{tr}

\newcommand{\AP}{\operatorname{Ap}}

\newcommand{\QF}{\operatorname{QF}}

\newcommand{\Ext}{\operatorname{Ext}}
\newcommand{\typ}{\operatorname{type}}
\newcommand{\Gr}{\operatorname{group}}

\def\0{{\mathbf 0}}

\def\c{{\mathbf c}}
\def\m{{\mathbf m}}
\def\1{{\mathbf 1}}
\def\a{{\mathbf a}}
\def\b{{\mathbf b}}
\def\z{{\mathbf z}}
\def\f{{\mathbf f}}

\def\h{{\mathbf h}}

\def\bs{{\mathbf s}}
\def\x{{\mathbf x}}

\def\v{{\mathbf v}}
\def\w{{\mathbf w}}

\def\fp{{\mathfrak p}}

\def\fm{{\mathfrak m}}

\def\fp{{\mathfrak p}}

\newtheorem{theorem}{Theorem}[section]
\newtheorem{corollary}[theorem]{Corollary}
\newtheorem{lemma}[theorem]{Lemma}
\newtheorem{proposition}[theorem]{Proposition}

\theoremstyle{definition}
\newtheorem{definition}[theorem]{Definition}
 
\newtheorem{example}[theorem]{Example}
\newtheorem{remark}[theorem]{Remark}

\numberwithin{equation}{section}

\title{On nearly Gorenstein affine semigroups}

\author{Raheleh Jafari}
\address{Mosaheb Institute of Mathematics, Kharazmi University, Tehran, Iran.}
\email{rjafari@khu.ac.ir}

\author{Francesco Strazzanti}
\address{Dipartimento di Matematica, Dipartimento di Eccellenza 2023-2027, Universit\`a di Genova, Via Dodecaneso 35, 16146 Genova, Italy.}
\email{francesco.strazzanti@gmail.com}

\author{Santiago Zarzuela Armengou}
\address{Departament de Matem\`{a}tiques i Inform\`{a}tica, Universitat de Barcelona, Gran via de les Corts Catalanes 585, 08007 Barcelona, Spain.}
\email{szarzuela@ub.edu}

\thanks{The authors received partial support from the grant PID2022-137283NB-C22 funded by 	MCIN/AEI/10.13039/501100011033. 
The second author is a member of the ``National
Group for Algebraic and Geometric Structures, and their Applications'' (GNSAGA - INdAM)}

\begin{document}
	
\begin{abstract}
We describe the canonical module of a simplicial affine semigroup ring $\KK[S]$ and its trace ideal. As a consequence, we characterize when $\KK[S]$ is nearly Gorenstein in terms of arithmetic properties of the semigroup $S$. Then, we find some bounds for the Cohen-Macaulay type of $\KK[S]$ when it is nearly Gorenstein. In particular, if it has codimension at most three, we prove that the Cohen-Macaulay type is at most three and this bound is sharp.  
\end{abstract}
	
	
	
\maketitle

\section*{Introduction}

The notion of nearly Gorenstein ring appeared several times in literature as a generalization of Gorenstein rings, even if this name was introduced later in \cite{HHS}, see for instance \cite{D,HV,SV}. Indeed, only in 2019, with the work of Herzog, Hibi, and Stamate \cite{HHS}, a systematic study of these rings has begun. Since then, many authors have studied this notion in several contexts like numerical semigroup rings \cite{HHS2,Moscariello-Strazzanti}, projective monomial curves \cite{Mi}, Ehrhart rings \cite{HKMM}, and quotient singularities \cite{CS} among others. The idea behind nearly Gorenstein rings relies on the fact that the trace ideal of the canonical module of a ring determines its non-Gorenstein locus. More precisely, let $(R,\mathfrak m)$ be a Cohen-Macaulay local ring (or a positively graded $\KK$-algebra) admitting canonical module $\omega_R$, and define the trace ideal of $\omega_R$ as the ideal
\[
\tr(\omega_R)= \sum_{\varphi \in \Hom_R(\omega_R,R)} \varphi(\omega_R).
\] 
Given a prime ideal $\mathfrak{p} \in {\rm Spec}(R)$, the ring $R_{\mathfrak{p}}$ is not Gorenstein if and only if $\tr(\omega_R) \subseteq \mathfrak{p}$, see \cite[Lemma 2.1]{HHS}. Therefore, $R$ is Gorenstein if and only if $\tr(\omega_R)=R$. The ring $R$ is said to be nearly Gorenstein if $\tr(\omega_R)$ contains the maximal (homogeneous) ideal $\mathfrak{m}$ of $R$, i.e., if $\tr(\omega_R)$ is equal to either $R$ or $\mathfrak{m}$. It is clear that a nearly Gorenstein ring is Gorenstein in the punctured spectrum, but the converse is not true; indeed, $R$ is Gorenstein on the punctured spectrum if and only if $\tr(\omega_R)$ is an $\mathfrak{m}$-primary ideal. 

In this paper we are interested in the nearly Gorenstein affine semigroup rings, and especially in their Cohen-Macaulay type, which we call only type for brevity.
An affine semigroup $S$ is an additive submonoid of $\mathbb{N}^d$ for some positive integer $d$. We assume that $S$ is simplicial and fully embedded in $\mathbb{N}^d$. Its associated affine semigroup ring is the subalgebra of the polynomial ring $\mathbb{K}[x_1,x_2, \dots, x_d]$ given by $\mathbb{K}[S]=\mathbb{K}[x_1^{a_1} x_2^{a_2} \cdots x_d^{a_d} : (a_1, \dots,a_d)\in S]$, where $\mathbb{K}$ is a field. 
The case $d=1$ corresponds (up to isomorphism) to the case of numerical semigroup rings, where the nearly Gorenstein property has been more studied. In this case, if $\mathbb{K}[S]$ has embedding dimension at most $3$, it is well known that the type of $\mathbb{K}[S]$ is at most $2$ \cite{H}. Moreover, if $\mathbb{K}[S]$ has embedding dimension $4$ and is nearly Gorenstein, then in \cite{Moscariello-Strazzanti} it is proved that its type is at most $3$. 

More generally, for any $d$, when the embedding dimension of $\mathbb{K}[S]$ is at most $d+2$, in \cite{JY} it is shown that the type of $\mathbb{K}[S]$ is at most $2$. In this paper we focus on nearly Gorenstein affine semigroup rings with embedding dimension $d+3$ and prove that their type is at most $3$. Notice that without assuming that $\mathbb{K}[S]$ is nearly Gorenstein, the type is not bounded, even when $d=1$. Moreover, even in the numerical semigroup case, it is not known if the type of a nearly Gorenstein ring with fixed embedding dimension is bounded, see \cite[Question 3.7]{S}. 

We also prove that, regardless of the embedding dimension, in the nearly Gorenstein case the type is at least $d$ if the ring is not Gorenstein. The starting point to prove these results is the descriptions of the canonical module and its trace ideal in terms of the maximal elements of an Ap\'{e}ry set of $S$ with respect to a suitable order.

The structure of the paper is the following. In the first section, after recalling some definitions and some results, we characterize the canonical module of an affine semigroup ring and its trace ideal, see Theorem \ref{canonical module} and Proposition \ref{qf:tr} respectively. Among the consequences of these results, there is also a characterization of affine semigroup rings that are nearly Gorenstein or Gorenstein on the punctured spectrum, see Corollary \ref{NG} and Proposition \ref{3eq}. In Section 2 we focus on nearly Gorenstein rings $\mathbb{K}[S]$ and their type. For instance, if $\mathbb{K}[S]$ is not Gorenstein, in Corollary \ref{dtyp} we prove that its type is at least $d$. Finally, Section 3 is devoted to prove that if $\mathbb{K}[S]$ has embedding dimension $d+3$, then its type is at most $3$, see Theorem \ref{thm}. Both these bounds are sharp, indeed every possible value among these bounds can be obtained.
	
Several computations are performed by using Macaulay2 \cite{M2}, the GAP system \cite{GAP} and, in particular, the NumericalSgps package \cite{DGM}.

\section{The Canonical module and its trace}

Throughout the paper, $S=\langle\a_1,\dots,\a_{d+r}\rangle\subseteq\NN^d$ will be a simplicial and fully embedded affine semigroup with extremal rays $\a_1,\dots,\a_d$. This means that the vectors  $\a_{1},\dots,\a_{d}$  are linearly independent and for each element $\a\in S$, we have $n\a\in\NN\a_{1}+\dots+\NN\a_{d}$, for some positive integer $n$. Equivalently
\begin{equation*}
S\subseteq\sum^d_{i=1}\QQ_{\geq0} \ \a_i.
\end{equation*}
Moreover, we assume that $\a_i$ is the smallest generator in its extremal ray for every $i=1,\dots,d$, i.e., if there is another generator $\a_j$ of $S$ for which $\a_j=q \a_i$ for some non-zero $q \in \mathbb{Q}$, then $q>1$. 

We write $\Gr(S)$ for the smallest group (up to isomorphism) that contains $S$, $\Gr(S)=\{\a-\b \mid \a, \b \in S\}.$ 
Given two sets $A,B\subseteq\NN^d$, we write $A+B$ for the set $\{\a+\b \ ; \ \a\in A , \b\in B\}$. If $A=\{\a\}$, we simply write $\a+B$, instead of $\{\a\}+B$.   
A subset $H\subseteq\Gr(S)$ is called $S$-ideal of $\Gr(S)$, if 
$S+H\subseteq H$.  
When $H\subseteq S$, we simply say that $H$ is an ideal of $S$.

Let $R=\KK[S]$ be the semigroup ring associated to $S$.  
A monomial in the semigroup ring $R$ is an element of
the form $\x^\a=x_1^{a_1}x_2^{a_2}\dots x_{d}^{a_{d}}$, where $\a=(a_1,a_2,\dots,a_{d})\in S$. An ideal $I\subseteq R$ is a monomial ideal if it is generated by monomials. For any subset $H$ of $S$, let  $\KK[H]$ denote the $\KK$-linear span of the
monomials $\x^\a$ with $\a\in H$. Then, $I$ is a monomial ideal if and only if $I=\KK[H]$ for some ideal $H$ of $S$, or equivalently, if $I$ is homogeneous with respect
to the tautological grading on $R$, which is defined by $\deg(\x^\a)=\a$.  Note that $\fm=\KK[M]$, where $M=S\setminus\{0\}$, is the unique monomial maximal ideal of $R$. 
Given a monomial fractional ideal  $I$ of $R$, let  $H$ be the set of exponents of monomials in $I$ and let $I^{-1}=\{x \in Q(R) : xI\subseteq R\}$, where $Q(R)$ is the field of fractions of $R$, and $H^{-1}=\{\z\in\Gr(S) \ ; \ \z+H\subseteq S \}$.
Then, $I=\KK[H]$ and $I^{-1}=\KK[H^{-1}]$.

The {\it Ap\'{e}ry set} of $S$ with respect to an element $\b\in S$ is defined as $\Ap(S,\b)=\{\a\in S \ ; \ \a-\b\notin S\}$. We will denote the zero vector of $\NN^d$ by $\0$. Since $S\subseteq\NN^d$, for $\b\neq \0$ we have $\0\in \Ap(S,\b)$. For a subset $A$, we set $\Ap(S,A)=\{\a\in S \ ; \ \a-\b\notin S, \text{ for all } \b\in A\}.$
Let $E=\{\a_1,\dots,\a_d\}$ and let $l_i$ be the smallest positive integer such that $l_i\a_{d+i}\in \sum^d_{j=1}\a_j$, for $i=1,\dots,r$. Then 
\[
\Ap(S,E)=\bigcap^d_{i=1}\Ap(S,\a_i)\subseteq\left\{\sum^r_{i=1}n_i\a_{d+i} \; ; \; 0\leq n_i<l_i\right\},
\]
is a finite set.  
We consider the natural partial ordering $\preceq_S$ on $S$ where, for all elements $\a$ and $\b$ in $\NN^d$, $\b\preceq_S\a$ if there is an element $\c\in S$ such that $\a=\b+\c$.  

Given $\b\in\max_{\preceq_S}\Ap(S,E)$, the element $\b-\sum^d_{i=1}\a_i$ is said to be a {\em quasi-Frobenius element} and the set of quasi-Frobenius elements of $S$ is denoted by $\QF(S)$, see \cite[Defenition~3.1]{JY}.
If $R$ is Cohen-Macaulay, then the number of quasi-Frobenius elements is equal to the Cohen-Macaulay type of $R$, \cite[Proposition~3.3]{JY}. We denote this number by $\typ(S)$ and refer to it as the {\it type} of the semigroup $S$.

When $R$ is Cohen-Macaulay, a finite graded $R$-module $C$ is a canonical module of $R$ if there exist homogeneous isomorphisms 
 \[
 \Ext^i_R(R/\fm,C)\cong\left\{ \begin{array}{ll} 0  & \text{ for } i\neq d, \\ R/\fm & \text{ for } i=d.   \end{array} \right.
 \]
It is unique up to isomorphism and $R$ is Gorenstein exactly when it is a canonical module of itself.
 
Let $F_i=({\sum^d_{\underset{j\neq i}{j=1}}}\QQ_+\a_j)\cap S$, and let \[G_i=\{\b\in\Gr(S) : \b+\a\in S \ \text{for some}\ \a\in F_i\}\]
for $i=1,\ldots,d$. Let $C_S=-(\cap_{i=1}^d C_i)$, where $C_i=\Gr(S)\setminus G_i$.
If $R$ is Cohen-Macaulay, then $\KK[C_S]$ is the canonical module of $R$ by \cite[Theorem~3.8]{GSW-76}. Let 
\[
\omega_S=\{-\f : \f\in\QF(S)\}+S\cong \{-\m : \m\in\max_{\preceq_S}\Ap(S,E)\}+S.
\]
In \cite[Definition~3.6]{W}, the authors define the graded canonical module of $\KK[S]$ to be $\KK[\omega_S]$ and show that in several aspects this definition is consistent with the  canonical module as generally defined. In the following we prove that this definition coincide with the abstract definition of a  canonical module. 
\begin{theorem}\label{canonical module}
If $R$ is Cohen-Macaulay, then $\KK[\omega_S]$ is a canonical module of $R$.
\end{theorem}
\begin{proof}
	Let $\m\in \max_{\preceq_S}\AP(S,E)$. First, we show that $\m-\sum^d_{i=1}\a_i\notin G_i$, for $i=1,\ldots,d$. Without loss of generality, we may assume that $i=1$. If $\m-\sum^d_{i=1}\a_i\in G_1$,  then $\m-\sum^d_{i=1}\a_i+\a\in S$, for some $\a=\sum^d_{i=2}\lambda_i\a_i\in S$, where $\lambda_i\in\QQ_+$, for $i=2,\ldots,d$. Let $l$ be a positive integer such that $l\lambda_i\in\NN$, for $i=2,\ldots,d$. Then $(l-1)\a\in S$, and so
	\[\m-\sum^d_{i=1}\a_i+l\a=\m-\a_1+\sum^d_{i=2}(l\lambda_i-1)\a_i\in S.\]
	Let $\bs=\m+\sum^d_{i=2}(l\lambda_i-1)\a_i=\a_1+\h$, for some $\h\in S$. Since $\m\in\AP(S,\a_1)$, we have $l\lambda_j\neq 1$, for some $j=2,\ldots,d$. Since $\m-\a_1\in \Gr(S)$, if $l\lambda_i-1$ and $l\lambda_j-1$ are positive for some $i\neq j$, then $\m-\a_1\in S$, by \cite[Theorem~2.6]{GSW-76}, a contradiction. Therefore,
	$\sum^d_{i=2}(l\lambda_i-1)\a_i=(l\lambda_j-1)\a_j$. Let $\alpha=l\lambda_j-1$. Then
	\[\bs=\m+\alpha \a_j=\a_1+\h,\]
	and \cite[Theorem~1.1]{RG-1998} implies that $\bs-\a_j-\a_1=\m-\a_1+(\alpha-1)\a_j\in S$. Applying \cite[Theorem~1.1]{RG-1998}, for $\alpha-1$ times, we get $\m-\a_1\in S$, a contradiction.
	
	Hence $\m-\sum^d_{j=1}\a_j\in\cap^d_{i=1}C_i$, and consequently, $\x^{-\f}$ belong to the canonical module of $R$ for every $\f \in \QF(S)$ by \cite[Theorem~3.8]{GSW-76}.

	Note that the embedding dimension of the canonical module of $R$ is equal to Cohen-Macaulay type of $R$  which is $|\QF(S)|$, by \cite[Proposition~3.3]{JY}. Therefore, it is enough to show that $\x^{-{\f}}$ cannot be generated by any other elements of $\KK[-\cap^d_{i=1}C_j]$, for $f \in \QF(S)$. Assume by contradiction that $\x^{-\f}=\x^{-\c}\x^{\bs}$, for some $\c\in\cap^d_{i=1}C_i$ and $s\in S$. Then $\c=\f+\bs$.
	Let $f=\m-\sum_{i=1}^d \a_i$. Because of the maximality of  $\m$, $\m+\bs-\a_k\in S$, for some $k$, $1\leq k\leq d$. It follows that $\c+\underset{\quad i\neq k}{\sum^d_{i=1}}\a_i=\f+\bs+\underset{\quad i\neq k}{\sum^d_{i=1}}\a_i=\m+\bs-\a_k\in S$, and so $\c\in G_i=\Gr(S)\setminus C_i$, a contradiction.
\end{proof}


We always denote the canonical module of $R$ by $\omega_R$.
For an $R$-module $N$, its {\em trace}, denoted by $\tr(N)$, is the sum of the ideals $\varphi(N)$ with $\varphi\in \Hom_{R}(N,R)$. Thus, $\tr(N)=\sum_{\varphi\in \Hom_{R}(N,R)}\varphi(N)$.
If $N_1\cong N_2$, then $\tr(N_1)=\tr(N_2)$, so while the canonical module $\omega_R$ is unique up to isomorphism, its trace is unique.
By \cite[Lemma 1.1]{HHS}, 
\begin{eqnarray*}
	\tr(\omega_R)&=&\omega_R \cdot (\omega_R)^{-1}\\
	&=& \KK[\omega_S] \cdot \KK[(\omega_S)^{-1}]\\
	&=&\KK[\omega_S+\omega_S^{-1}].
\end{eqnarray*}
We define the trace of $S$, to be the ideal $\tr(S)\subseteq S$ consisting of the exponents of monomials in $\tr(\omega_R)$. Then
\[
\tr(S)=\omega_S+\omega_S^{-1}.
\]

If $\KK[S]$ is Cohen-Macaulay, it is well known that $\tr(S)=S$ if and only if $\KK[S]$ is a Gorenstein ring. Moreover, the following result holds:

\begin{proposition}\cite[Lemma~2.1]{HHS}\label{tr}
Let $A$ be a Cohen-Macaulay finitely generated positively graded $\KK$-algebra, where $\KK$ is a field. For a prime ideal $\fp$ of $A$, the ring $A_\fp$ is not a Gorenstein ring if and only if $\tr(\omega_A)\subseteq\fp$.	
\end{proposition}

The following notion arises from these observations.

\begin{definition}
Assume that $A$ is a Cohen-Macaulay finitely generated positively graded $\KK$-algebra, where $\KK$ is a field. Let $\fm$ be the graded maximal ideal of $A$ and $\omega_A$ a canonical module of $A$. Then, $A$ is {\it nearly Gorenstein} if $\fm \subseteq \tr(\omega_{A})$.
\end{definition}

\begin{definition}
Assume that $\KK[S]$ is Cohen-Macaulay. We say that $S$ is {\it nearly Gorenstein} if $M \subseteq\tr(S)$, where $M=\langle \a_1, \a_2, \dots, \a_{d+r}\rangle$ is the maximal ideal of $S$. This is equivalent to require that $\KK[S]$ is nearly Gorenstein. 
\end{definition}

We can give a description of the elements of $\tr(S)$, which will yield a characterization of the nearly Gorenstein property.

\begin{proposition} \label{qf:tr}
Let $\max_{\preceq_S}\Ap(S,E)=\{\m_1,\dots,\m_t\}$. 
Then 
\[
\tr(S)=\{\b\in S \ ; \text{ there exists } 1\leq i\leq t  \text{ s.t. } \b+\m_i-\m_j\in S ~\text{for all}~ j=1, \dots, t\}.	
\]
\end{proposition}

\begin{proof}
Let $\Gamma$ be the set on the right side. Given $\b\in \Gamma$, there is $1\leq i\leq t$ such that $\b+\m_i-\m_j\in S$ for all $j=1, \dots, t$. 
	Thus, $\b+\m_i-\sum_{i=1}^d \a_d\in\omega_S^{-1}$ and $\b=-(\m_i-\sum_{i=1}^d \a_d)+\m_i-\sum_{i=1}^d \a_d+\b\in\tr(S)$.
	
	Now let $\a\in\tr(S)$, so there exist $\b_1\in\omega_S$ and $\b_2\in\omega_S^{-1}$ such that $\a=\b_1+\b_2.$
	Since $\b_1\in \omega_S$, we have $\b_1=-(\m_i-\sum_{i=1}^d \a_d)+\s$ for some $\s\in S$ and $1\leq i \leq t$.
	Since $\b_2\in \omega_S^{-1}$, we have $\b_2-(\m_j-\sum_{i=1}^d \a_d)\in S$ for every $j=1,\dots,t$. Therefore
	\begin{center}
		$\a+\m_i-\m_j=\b_1+\b_2+\m_i-\m_j=\s+\b_2-(\m_j-\sum_{i=1}^d \a_d)\in S.$
	\end{center}
	Hence, $\tr(S)\subseteq \Gamma$.	
\end{proof}

The following corollary is a generalization of \cite[Proposition 1.1]{Moscariello-Strazzanti}, which was proved in the case $d=1$.

\begin{corollary} \label{NG}
Assume that $\KK[S]$ is Cohen-Macaulay and let $\max_{\preceq_S}\Ap(S,E)=\{\m_1,\dots,\m_t\}$. 
Then, $S$ is nearly Gorenstein if and only if for each $i=1, \dots, d+r$ there exists $\m \in \max_{\preceq_S}\Ap(S,E)$ such that $\a_i+\m-\m_j\in S$ for each $j=1, \dots, t$.
\end{corollary}

\begin{remark}\label{2max}
If $S$ has type $2$ and $\max_{\preceq_S}\Ap(S,E)=\{\m_1,\m_2\}$, then 
\[
\tr(S)=\{\b\in S \ ; \ \b+\m_1-\m_2\in S \ or \ \b+\m_2-\m_1~\in S\}.
\]
In particular, we always have $\{\m_1,\m_2\}\subseteq\tr(S)$.
\end{remark}


In the previous results we have assumed that $\KK[S]$ is Cohen-Macaulay. However, this property can be read off from the semigroup $S$, as the next result shows. We denote with $\Gr(\a_1, \dots, \a_d)$ the group generated by $\a_1, \dots, \a_d$, i.e., $\Gr(\a_1, \dots, \a_d)=\Gr(\langle \a_1, \dots, \a_d\rangle)$.

\begin{proposition}\cite[Theorem 1.5 and Corollary 1.6]{RG-1998}\label{CM}
	The following statements are equivalent:
	\begin{enumerate}
		\item $\KK[S]$ is Cohen-Macaulay;
		\item For all $\w_1,\w_2\in\Ap(S,E)$, if $\w_1-\w_2\in\Gr(\a_1,\dots,\a_d)$, then $\w_1=\w_2$; 
		\item For all $\w_1,\w_2\in\Ap(S,E)$, if $\w_1+\sum^d_{i=1}l_i\a_i=\w_2+\sum^d_{i=1}l'_i\a_i$ with $l_i, l_i' \in \mathbb{N}$, then $\w_1=\w_2$ and $l_i=l'_i$ for $i=1,\dots,d$.
	\end{enumerate}	
\end{proposition}

We now start to explore some consequences of Proposition \ref{qf:tr}.

\begin{proposition}\label{hplane}
Assume that $\KK[S]$ is Cohen-Macaulay.	If there is an hyperplane containing all the  generators $\a_{d+1},\dots,\a_{d+r}$ and the origin of coordinates, then the following statements are equivalent:
\begin{enumerate}
	\item $\a_1,\dots,\a_d\in\tr(S)$; 
	\item $\KK[S]$ is Gorenstein;
	\item $\tr(S)=S$;
	\item $S$ is nearly Gorenstein. 
\end{enumerate}
\end{proposition}
\begin{proof}
Let $\DD$ denote the $(d-1)$-dimensional space that contains  $\a_{d+1},\dots,\a_{d+r}$ and the origin of coordinates. 
Note that all the elements of $\Ap(S,E)$ are on $\DD$ and at least one extremal ray $\a_i$, does not belong to this hyperplane. By Proposition \ref{qf:tr}, there exists $T(\a_i)\in \max_{\preceq_S}\Ap(S,E)$ such that $\a_i+T(\a_i)-\m \in S$ for every $\m \in \max_{\preceq_S}\Ap(S,E)$. 

(1)$\implies$(2)
Assume that  $\typ(S)>1$ and consider $\m\in\max_{\preceq_S}\Ap(S,E)\setminus\{T(\a_i)\}$. As $\a_i$ is not on $\DD$, we have  $T(\a_i)+\a_i-\m\notin\Ap(S,E)$ and then there exist $\a_j$ with $1 \leq j \leq d$ and $s \in S$ such that $\a_i+T(\a_i)-\m=\a_j+s$. Writing $\m+s=\sum_{k=1}^{d}l_k \a_k+\n$ for some $l_k \in \mathbb{N}$ and $\n \in \Ap(S,E)$, we get $\a_i+T(\a_i)=\a_j+\sum_{k=1}^{d}l_k \a_k+\n$. By Proposition \ref{CM}, $\sum_{k=1}^{d}l_k \a_k=0$ and $i=j$, so that $\m=\n=T(\a_i)$, which is a contradiction. 

(2)$\implies$(3),  (3)$\implies$(4), and (4)$\implies$(1) are clear.
\end{proof}

The same proof of the previous proposition gives the following result that we will repeatedly use in the next sections. 

\begin{corollary}\label{line}
Let $\KK[S]$ be Cohen-Macaulay. If all the generators $\a_{d+1},\dots,\a_{d+r}$ are on the same line passing through the origin of coordinates, then the following statements are equivalent:
	\begin{enumerate}
		\item $\a_i,\a_j\in\tr(S)$ for some $1\leq i<j\leq d$;
		\item $\KK[S]$ is Gorenstein;
		\item $\tr(S)=S$;
		\item $S$ is nearly Gorenstein. 
	\end{enumerate}
\end{corollary}

Given $\w \in \Ap(S,E)$, we say that $\w$ has a {\it unique expression} if it can be written as $\w=\sum_{i=1}^{r} l_i\a_{d+i}$ with $l_i \in \mathbb{N}$ in a unique way.

\begin{corollary}\label{cuniq}
	Assume that $\KK[S]$ is Cohen-Macaulay and that $S$ has embedding dimension $2d$. If $E\subseteq\tr(S)\subseteq M$, then the non-extremal generators $\a_{d+1},\dots,\a_{2d}$ are linearly independent. In particular, each element of $\Ap(S,E)$ has a unique expression.
\end{corollary}
\begin{proof}
	Note that the cone generated by $\a_{d+1},\dots,\a_{2d}$ is a $d$ dimensional space by Proposition~\ref{hplane}. This means that there is no relation between them, in particular all the elements in $\Ap(S,E)$ have unique expressions. 
\end{proof}

We end this section by characterizing when $\KK[S]$ is Gorenstein on the punctured spectrum. This result generalizes \cite[Corollary~3.4]{JTY}.

\begin{proposition}\label{3eq}
Assume that $\KK[S]$ is Cohen-Macaulay. The following conditions are equivalent:
	\begin{enumerate}
		\item $\KK[S]$ is Gorenstein on the punctured spectrum;
		\item $\lambda_1\a_1,\dots,\lambda_d\a_d\in\tr(S)$ for some $\lambda_1,\dots,\lambda_d\in\NN$;
		\item There exist $\lambda_1,\dots,\lambda_d\in\NN$ and $1\leq r_1,\dots,r_d\leq t$ such that 
		$\lambda_i\a_i+\m_{r_i}-\m_j\in S$ for all $1\leq j\leq t$;
	\item For any $1\leq i\leq d$, there exist $\w_i\in\Ap(S,E)$ and $\lambda_i\in\NN$ such that $$\w\preceq\w_i+\lambda_i\a_i,$$ for all $\w\in\Ap(S,E).$
\end{enumerate}
	\end{proposition}
\begin{proof}
	(1)$\implies$(2) It  follows by Proposition~\ref{tr} and the fact that $\x^{\a_1},\dots,\x^{\a_d}$ provide a system of parameters for $\KK[S]$.

	(2)$\implies$(3) It follows by Proposition~\ref{qf:tr}.

	(3)$\implies$(4) This is clear.
	
	(4)$\implies$(1) Let $\w_i\preceq_S\m_{r_i}$ for $i=1,\dots,d$. As $\m_{r_i}\preceq\w_i+\lambda_i\a_i$, we get $\w_i=\m_{r_i}$. Therefore, $\lambda_i\a_i\in\tr(S)$ by Proposition~\ref{qf:tr}. Since $\x^{\a_1},\dots,\x^{\a_d}$ provide a system of parameters, it follows that $\tr(\omega_R)$ is an $\fm$-primary ideal. Now, Proposition~\ref{tr} implies (1).
\end{proof}

By Proposition \ref{tr}, it is clear that a nearly Gorenstein ring is Gorenstein on the punctured spectrum. However, the converse is false, even for affine semigroup rings.

\begin{example} \label{not NG}
Let $S$ be the semigroup generated by $\a_1=(6,0), \a_2=(0,6), \a_3=(2,1), \a_4=(1,2)$. In this case $\Ap(S,E)=\{\0, \a_3, 2\a_3, 3\a_3, \a_4, 2\a_4, 3\a_4, \m_1=3\a_3+\a_4, \m_2=\a_3+3\a_4\}$ and by Proposition \ref{CM} one can see that $\KK[S]$ is Cohen-Macaulay. Moreover, we have that $\a_1+\m_2-\m_1=2\a_3$ and $\a_2+\m_1-\m_2=2\a_4$ belong to $S$, whereas a simple computation shows that $S\setminus\tr(S)=\{0,\a_3,\a_4,\a_3+\a_4\}$. Therefore, $\KK[S]$ is Gorenstein on the punctured spectrum, but not nearly Gorenstein.
\end{example}

\section{Nearly Gorenstein affine semigroups}

In this section we focus on nearly Gorenstein affine semigroups and their type. As in the previous section, $S$ will always be a simplicial affine semigroup fully embedded in $\mathbb{N}^d$ and minimally generated by $\a_1, \a_2, \dots, \a_{d+r}$, whose extremal rays are $\a_1, \dots, \a_d$. Recall that $S$ is nearly Gorenstein if $\a_i \in \tr(S)$ for every $i=1, \dots, d+r$. Moreover, by Proposition \ref{qf:tr}, $\a_i \in \tr(S)$ exactly when there exists $\m\in\max_{\preceq_S}\Ap(S,E)$ such that  $\m+\a_i-\n\in S$ for all $\n\in\max_{\preceq_S}\Ap(S,E)$. In the following lemma we show that if this happens for an extremal ray, then $\m$ is unique.

\begin{lemma}\label{uni-m}
Let $\a_i\in\tr(S)$ for some $1\leq i\leq d$. Then, there exists a unique $\m\in\max_{\preceq_S}\Ap(S,E)$ such that  $\m+\a_i-\n\in S$ for all $\n\in\max_{\preceq_S}\Ap(S,E)$.
\end{lemma}
\begin{proof}
	Assume contrary that there are two maximal elements $\m,\n$ with this property. Then 
	\[
	\m+\a_i=\n+\w  \ , \ \n+\a_i=\m+\w'
	\]
	for some $\w,\w'\in S$. Thus
	\[
	\m+2\a_i=\n+\a_i+\w=\m+\w+\w',
	\]
	which implies that $2\a_i=\w+\w'$. As $\a_i$ is on an extremal ray of the cone generated by $S$, it follows that $\w=l\a_i$ and $\w'=l'\a_i$ for some $l,l'\in\QQ$ with $l+l'=2$. Thus, either $l\leq 1$ or $l'\leq 1$. Since $\a_i$ is the smallest element on this extremal ray, we have $\w=\a_i$ and $\w'=\a_i$ and consequently $\m=\n$. 
\end{proof}

\begin{example}
In the previous lemma we need that $\a_i$ is the smallest element in its extremal ray, otherwise the statement is false.
For instance, let $S$ be the semigroup generated by $\a_1=(2,0)$, $\a_2=(0,2)$, $\a_3=(0,3)$, $\a_4=(1,1)$, $\a_5=(1,2)$. We have $\max_{\preceq_S}\Ap(S,E)=\{\m_1=(1,1),\m_2=(1,2),\m_3=(0,3)\}$ and
\begin{align*}
&\m_3+\a_1-\m_1\in S \ &\m_3+\a_1-\m_2\in S \ \ \ \ \ \ \ &\m_1+\a_1-\m_2\notin S \ &\m_2+\a_1-\m_3\notin S, \\
&\m_2+\a_2-\m_1\in S \ &\m_2+\a_2-\m_3\in S \ \ \ \ \ \ \ &\m_1+\a_2-\m_2\notin S \ &\m_3+\a_2-\m_1\notin S.
\end{align*}
Therefore, as we proved in the previous lemma, for $i=1,2$ there is only one $\m\in\max_{\preceq_S}\Ap(S,E)$ for which $\m+\a_i-\n\in S$ for all $\n\in\max_{\preceq_S}\Ap(S,E)$. However, this is not true for $\a_3$, even though it is in the same line of $\a_2$. Indeed, both $\m_1$ and $\m_2$ work for $\a_3$.

Using Proposition \ref{qf:tr} and Proposition \ref{CM}, it is not difficult to see that $S$ is nearly Gorenstein. Moreover, also for $\a_4$ and $\a_5$ there are two maximal elements of $\Ap(S,E)$ for which the equalities above hold.
\end{example}

Throughout the paper, $T(\a_i)$ will denote the unique maximal Ap\'ery element assigned to $\a_i$ in Lemma~\ref{uni-m} for $i=1,\dots,d$.


\begin{proposition}\label{inj-T}
Let $\KK[S]$ be Cohen-Macaulay with $\typ(S)\geq2$. If   $1\leq r<s\leq d$, then $T(\a_r)\neq T(\a_s)$.
\end{proposition}
\begin{proof}
Assume on the contrary that $T(\a_r)=T(\a_s)$ and set $\m=T(\a_r)=T(\a_s)$. Since $\typ(S)\geq2$, we may choose $\n\in\max_{\preceq_S}\Ap(S,E)\setminus\{\m\}$.
	Then 
$\b_r=\m+\a_r-\n$ and $\b_s=\m+\a_s-\n$ belong to $S$.	
Let $\b_r=\w_r+\sum^d_{i=1}l_i\a_i$ and  $\b_s=\w_s+\sum^d_{i=1}l'_i\a_i$, where $\w_r,\w_s\in\Ap(S,E)$. Note that 
\[
\b_r-\b_s=\w_r-\w_s+\sum^d_{i=1}(l_i-l'_i)\a_i=\a_r-\a_s.
\]
Since $S$ is Cohen-Macaulay, Proposition~\ref{CM} implies that $\w_r=\w_s$ and $\sum^d_{i=1}l_i\a_i+\a_s=\sum^d_{i=1}l'_i\a_i+\a_r$. Since $\a_1, \dots \a_d$ are linearly independent, $l_r=l'_r+1>0$. As $\m-\n+\a_r=\w_r+\sum^d_{i=1}l_i\a_i$, we get $\m-\n=\w_r+\sum_{i\neq r}l_i\a_i+(l_r-1)\a_r\in S$, which is a contradiction. 
\end{proof}

As a consequence, if $S$ is nearly Gorenstein and $\typ(S) \geq 2$, then there exists at least $d$ different maximal elements in $\Ap(S,E)$.

\begin{corollary}\label{dtyp}
If $S$ is nearly Gorenstein but not Gorenstein, then $\typ(S)\geq d$.
\end{corollary}

\begin{example} \label{d=t=2}
The bound in the previous corollary is sharp. For instance, let $S$ be the semigroup generated by $\a_1=(0,3)$, $\a_2=(3,1)$, $\a_3=(1,2)$, $\a_4=(2,2)$, $\a_5=(3,3)$. In this case $\typ(S)=2=d$ because $\max_{\preceq_S}\Ap(S,E)=\{(5,7),(6,6)\}$. By Proposition \ref{CM}, $\KK[S]$ is Cohen-Macaulay. Moreover, by applying Remark~\ref{2max}, it is straightforward to see that all the generators are in $\tr(S)$, and then $S$ is nearly Gorenstein. 
\end{example}

We have already seen that being $\KK[S]$ Gorenstein on the punctured spectrum does not imply that $\KK[S]$ is nearly Gorenstein. Notice that this happens even if $\typ(S)=d$, indeed in Example \ref{not NG} we have $\typ(S)=d=2$.
 

We can also say something about the embedding dimension of a nearly Gorenstein semigroup, but we first need a technical lemma. 

\begin{lemma}\label{lll}
Let $\KK[S]$ be Cohen-Macaulay,  $\m\in\max_{\preceq_S}\Ap(S,E)$ and $\v\in\Ap(S,E)$ such that $\v+\a_i=\m+\sum^r_{t=1}l_t\a_{d+t}$ for some  $1\leq i\leq d$ and $l_1,\dots,l_r\in\NN$. If $l_j\neq 0$, then $\m+\a_{d+j}=\w+\a_i$, where $\w\in\Ap(S,E)$ and $\w-\a_{d+j}\notin S$.  In particular, $\v-(\sum^{j-1}_{t=1}l_t\a_{d+t}+(l_j-1)\a_{d+j}+\sum^{r}_{t=j+1}l_t\a_{d+t})\in S$.
\end{lemma}
\begin{proof}
	As $\m\in\max_{\preceq_S}\Ap(S,E)$,  $\m+\a_{d+j}\notin\Ap(S,E)$ and we get 
	$\m+\a_{d+j}=\w+\sum^d_{t=1}h_t\a_t$ for some $\w\in\Ap(S,E)$. 
	Then 
	$$\v+\a_i=\w+\sum^d_{t=1}h_t\a_t+\underset{t\neq j}{\sum^r_{t=1}}l_t\a_{d+t}+(l_j-1)\a_{d+j}=\sum^d_{t=1}h_t\a_t+\w'+\sum^d_{t=1}h'_t\a_t,$$
	 for some $\w'\in\Ap(S,E)$ and $h_t,h'_t\in\NN$. Since $\KK[S]$ is Cohen-Macaulay and $\v,\w'\in\Ap(S,E)$, by Proposition~\ref{CM} we have that $\sum^d_{t=1}(h_t+h'_t)\a_t=\a_i$, which implies $\sum^d_{t=1}h_t\a_t=\a_i$. Therefore,  $\m+\a_{d+j}=\w+\a_i$.  Finally,  $\w-\a_{d+j}=\m-\a_i\notin S$ and $\v-(\sum^{j-1}_{t=1}l_t\a_{d+t}+(l_j-1)\a_{d+j}+\sum^{r}_{t=j+1}l_t\a_{d+t})=\m+\a_{d+j}-\a_i=\w\in S$. 
\end{proof}

\begin{proposition}\label{typ-emd}
The following statements hold when $S$ is  nearly Gorenstein.
\begin{enumerate}
	\item If $\typ(S)=d$, then $\edim(S)\geq 2d-1$.
	\item If $\typ(S)>d$, then $\edim(S)\geq 2d$.
\end{enumerate}
\end{proposition}
\begin{proof}
	Let $\m_{i}=T(\a_i)$ for $i=1,\dots,d$. Take $j \in \{1\dots d\}$ and  $\m\in\max_{\preceq_S}\Ap(S,E)\setminus \{\m_j\}$. Then
	\[
	\m_j+\a_j=\m+\sum^{r}_{i=1}l_{j,i}\, \a_{d+i}.
	\]
If $l_{j,i}\neq 0$, then $\m+\a_{d+i}=\w+\a_j$ for some $\w\in\Ap(S,E)$,  by Lemma~\ref{lll}. If $\typ(S)>d$ and we choose $\m \in \max_{\preceq_S}\Ap(S,E)\setminus \{\m_1, \dots, \m_d\}$, by Proposition \ref{CM} it is not possible to have $\w+\a_j=w'+\a_{j'}$ with $\w,\w' \in \Ap(S,E)$ and $j \neq j'$. Therefore, the number of elements in the set $\{i ; l_{j,i}>0  \text{ for some } j\}$ is $d$ and this implies that there are at least other $d$ minimal generators. The case in which $\typ(S)=d-1$ is similar. 
\end{proof}

By Proposition~\ref{typ-emd} and \cite[Theorem~3.5]{JY} we immediately get the following corollary.

\begin{corollary}
	If $S$ is nearly Gorenstein with $d\leq 3$, then $\edim(S)\geq 2d$.
\end{corollary}

When there is at least one $T(\a_i)$ having a unique expression, it is also possible to give an upper bound for $\typ(S)$. This will follow from the next lemma.

\begin{lemma}\label{uni-ex-li}
		Let $\KK[S]$ be Cohen-Macaulay and assume that $\m\in S$  has a unique expression. 
If $\m+\a_t=\n+\sum^r_{i=1}l_i\a_{d+i}$ for some $1\leq t\leq d$, $\n\in\max_{\preceq_S}\Ap(S,E)$ and $l_1,\dots,l_r\in\NN$, then there is only one $i$, $1\leq i\leq r$, for which $l_i\neq 0$.		
\end{lemma}
\begin{proof}
	Let $\m$  have the unique expression $\sum^r_{j=1}\lambda_j\a_{d+j}$. Then $\lambda_j=\max\{l \ ; \  \m -l\a_{d+j}\in S\}$.  
			If $l_j\neq0$, then $\n+\a_{d+j}=\w+\a_i$ for some $\w\in\Ap(S,E)$, by Lemma~\ref{lll}, which implies that $\m=\w+(l_j-1)\a_{d+j}+l_1\a_{d+1}+\dots+\widehat{l_j\a_{d+j}}+\dots+l_r\a_{d+r}$. Therefore, $l_j\leq\lambda_j$ for $j\neq i$. Now, if there  are $1\leq t\neq k\leq r$ with $l_t>0, l_k>0$, then $l_j\leq\lambda_j$ for all $j=1,\dots,r$ which implies $\n-\a_t\in S$, a contradiction. Thus, there is only one $1\leq j\leq r$, with $l_j\neq0$. 
\end{proof}

As a consequence we get the following result.

\begin{proposition}\label{uni-ex-r}
If $S$ is nearly Gorenstein and $T(\a_i)$ has a unique expression for some $1\leq i\leq d$, then $\typ(S)\leq r+1.$
\end{proposition}

We conclude this section collecting some properties that we will repeatedly use in the next section in the particular case $r=3$. 

\begin{lemma}\label{lem-T}
Assume that $S$ is nearly Gorenstein and let $i \in \{1, \dots, d\}$. For every $\m\in\max_{\preceq_S}\Ap(S,E)\setminus\{T(\a_i)\}$ there exist $\lambda^\m_{i,1},\dots,\lambda^\m_{i,r}\in\NN$ such that
	\[
	T(\a_i)+\a_i=\m+\sum^{r}_{s=1}\lambda^\m_{i,s}\a_{d+s}.
	\]
	\begin{enumerate}
		\item If $\lambda^\m_{i,s}>0$, for some $1\leq s\leq r$, then $\lambda^\m_{j,s}=0$ for  $j\in\{1,\dots,d\}\setminus\{i\}$ with $\m\neq T(\a_j)$.
		\item If there is only one $s$ with $\lambda^\m_{i,s}\neq0$, then 
		$\lambda^\m_{i,s}-1=\max\{l ; T(\a_i)-l\a_{d+s}\in S\}$. 
		\item Let  $\m,\n\in\max_{\preceq_S}\Ap(S,E)\setminus\{T(\a_i),T(\a_j)\}$ for some $1\leq i,j\leq d$. Then 
		\[
		\sum^{r}_{s=1}\lambda^\m_{i,s}\a_{d+s}+\sum^{r}_{s=1}\lambda^\n_{j,s}\a_{d+s}=\sum^{r}_{s=1}\lambda^\n_{i,s}\a_{d+s}+\sum^{r}_{s=1}\lambda^\m_{j,s}\a_{d+s}.
		\]
		\item Let $\m=T(\a_i)$ and $\n=T(\a_j)$ for  some $1\leq i,j\leq d$. Then
		\[
		\a_i+\a_j=\sum^{r}_{s=1}\lambda^\n_{i,s}\a_{d+s}+\sum^{r}_{s=1}\lambda^\m_{j,s}\a_{d+s}.
		\]
		\item Let $1\leq i\neq j\leq d$, $\n=T(\a_j)$ and $\m\in\max_{\preceq_S}\Ap(S,E)\setminus\{T(\a_i),T(\a_j)\}$. Then
		\[
		\sum^{r}_{s=1}\lambda^\n_{i,s}\a_{d+s}+\sum^{r}_{s=1}\lambda^\m_{j,s}\a_{d+s}=\sum^{r}_{s=1}\lambda^\m_{i,s}\a_{d+s}+\a_j.
		\]
	\end{enumerate}
\end{lemma}

\begin{proof}
	(1) By Lemma~\ref{lll},  
	$\m+\a_{d+s}=\w+\a_i$ for some $\w\in\Ap(S,E)$.  	
	If $\lambda^\m_{j,s}>0$, for some $j$, then using again Lemma~\ref{lll}, we get 	$\m+\a_{d+s}=\v+\a_j$ for some  $\v\in\Ap(S,E)$. Therefore, $i=j$ by Proposition \ref{CM}.
	
	(2) We have $T(\a_i)+\a_i=\m + \lambda_{i,s}^{\m}\a_{d+s}$ with $\lambda_{i,s}^{\m}>0$. As in (1), it follows that $\m + \a_{d+s}=\w+\a_i$ with $\w \in \Ap(S,E)$, and then $T(\a_i)=\w+(\lambda_{i,s}^{\m}-1)\a_{d+s}$. This means that $T(\a_i)-(\lambda_{i,s}^{\m}-1)\a_{d+s} \in S$. On the other hand,  $T(\a_i)-\lambda_{i,s}^{\m}\a_{d+s}$ is not in $S$ because $T(\a_1)-\lambda_{i,s}^{\m}\a_{d+s}=\m-\a_i$ and $\m \in \Ap(S,E)$.
	
	(3), (4), and (5) are easy computations.
\end{proof}

Let $\mathfrak{m}$ be the monomial maximal ideal of $\KK[S]$. Recall that an ideal $J$ of $\KK[S]$ is said to be a reduction of $\mathfrak{m}$ if $\mathfrak{m}^{n+1}=\mathfrak{m}^nJ$ for some $n \in \NN$, and a reduction is called minimal if there are no other reductions contained in it. The case in which $\mathfrak{m}$ admits a monomial minimal reduction has been studied in \cite{DJS}, where several properties of $\KK[S]$ and its associated graded ring have been characterized. As a consequence of the previous lemma, in this case a nearly Gorenstein ring is also Gorenstein.

\begin{corollary}
Let $\KK[S]$ be nearly Gorenstein. If the monomial maximal ideal of $\KK[S]$ has a monomial minimal reduction, then $\KK[S]$ is Gorenstein.
\end{corollary}
\begin{proof}
It follows by Lemma~\ref{lem-T}(4) and \cite[Lemma 2.4 and Theorem 3.2]{DJS}.
\end{proof}

\section{Nearly Gorenstein semigroups with codimension three}

Assume that $\KK[S]$ is Cohen-Macaulay. If the codimension of $S$ is at most two, then its type is either one or two, see \cite[Theorem 3.5]{JY}. On the other hand, if the codimension is three, the type can be arbitrarily large as showed in \cite[Example 3.8]{JY}. However, when $d=1$ and $\KK[S]$ is nearly Gorenstein, in \cite{Moscariello-Strazzanti} it has been proved that the type is at most three (see also \cite{M} for the almost Gorenstein case). 
In this section we focus on nearly Gorenstein simplicial affine semigroups with codimension three, which means that $S$ is minimally generated by the extremal rays $\a_1, \dots, \a_d$ and three more generators $\a_{d+1}, \a_{d+2}, \a_{d+3}$. In particular, our goal is to prove that the type of such a semigroup is at most three. We start with a general lemma about elements of $\NN^d$.

\begin{lemma}\label{con}
Let $\b_1,\b_2,\b_3\in\NN^d$. If there exist positive integers  $\lambda_i,\lambda_j,\mu_i,\mu_k,\gamma_j,\gamma_k$  and non-negative integers $l_k,l_k',l_j,l'_j,l_i.l'_i$ such that 
	\begin{eqnarray*}
		\lambda_i\b_i+l_k\b_k=\lambda_j\b_j+l'_k\b_k\\
		\mu_i\b_i+l_j\b_j=\mu_k\b_k+l'_j\b_j\\
		\gamma_j\b_j+l_i\b_i=\gamma_k\b_k+l'_i\b_i,
	\end{eqnarray*}
	then $\b_i,\b_j,\b_k$ are on the same line passing   through the origin of coordinates.
\end{lemma}
\begin{proof}
	If  $l_k\geq l'_k$, then the first equation implies that $\b_j$ belongs to the cone   generated by $\b_i,\b_k$. 
	By the second equation, we have either $\b_i$ belongs to the cone generated by $\b_j,\b_k$, or  $\b_k$ belongs to the cone generated by $\b_i,\b_j$.
	In both cases it follows that  $\b_i,\b_j,\b_k$ are on the same line passing   through the origin.
	
	Now, assume that $l_k<l'_k$. Then  the first equation implies that $\b_i$ belongs to the cone   generated by $\b_j,\b_k$. 
	By the third equation, we have either $\b_j$ belongs to the cone generated by $\b_i,\b_k$, or  $\b_k$ belongs to the cone generated by $\b_i,\b_j$.
 In both cases it follows that  $\b_i,\b_j,\b_k$ are on the same line passing through the origin.
\end{proof}

By Corollary \ref{NG}, for any $\m \in \Ap(S,E)$ and $1 \leq i \leq d$ we have $T(\a_i)+\a_i=\m+\sum^3_{s=1}\lambda_s\a_{d+s}$ with $\lambda_s \in \NN$ for $s=1,2,3$. Having a deep understanding of these possible writings will be crucial in order to count the possible $\m \in \Ap(S,E)$. We start by showing that it is not possible that $\lambda_s$ is positive for all $s$.

\begin{corollary}\label{j=0}
	Let $S$ be nearly Gorenstein  of embedding dimension $d+3$ that  is not Gorenstein and let $1\leq i\leq d$. If $T(\a_i)+\a_i=\m+\sum^3_{s=1}\lambda_s\a_{d+s}$, for some $\m\in\max_{\preceq_S}\Ap(S,E)$, then there exists   $1\leq s\leq 3$, with $\lambda_s=0$.
\end{corollary}
\begin{proof}
	Assume on the contrary that $\lambda_s>0$ for $s=1,2,3$.  Then
			 $\m+\a_{d+s}=\w_s+\a_i$, for some $\w_s\in\Ap(S,E)$, by Lemma~\ref{lll}. Therefore, 
	\begin{eqnarray*}
		T(\a_i)&=&\w_1+(\lambda_1-1)\a_{d+1}+\lambda_2\a_{d+2}+\lambda_3\a_{d+3}\\
		&=&\w_2+\lambda_1\a_{d+1}+(\lambda_2-1)\a_{d+2}+\lambda_3\a_{d+3}\\&=&\w_3+\lambda_1\a_{d+1}+\lambda_2\a_{d+2}+(\lambda_3-1)\a_{d+3}.
	\end{eqnarray*}
	Note that $\w_s-\a_{d+s}\notin S$. Let $\w_1=l_2\a_{d+2}+l_3\a_{d+3}$, $\w_2=h_1\a_{d+1}+h_3\a_{d+3}$ and $\w_3=\mu_1\a_{d+1}+\mu_2\a_{d+2}$. Then 
	\begin{eqnarray*}
		(l_2+1)\a_{d+2}+(l_3+\lambda_3)\a_{d+3}&=&(h_1+1)\a_{d+1}+(h_3+\lambda_3)\a_{d+3}\\
		(l_2+\lambda_2)\a_{d+2}+(l_3+1)\a_{d+3}&=&(\mu_1+1)\a_{d+1}+(\mu_2+\lambda_2)\a_{d+2}\\
		(h_1+\lambda_1)\a_{d+1}+(h_3+1)\a_{d+3}&=&(\mu_1+\lambda_1)\a_{d+1}+(\mu_2+1)\a_{d+2}.
	\end{eqnarray*}	
	Now, 	Lemma~\ref{con}  implies that $\a_{d+1},\a_{d+2},\a_{d+3}$ are on the same line passing   through the origin of coordinates, a contradiction by Corollary~\ref{line}.
\end{proof}
  
If there are exactly two positive integers among $\lambda_1, \lambda_2$, and $\lambda_3$, we are not so lucky. However, in the next lemma we prove that for every choice of a couple of indices there is only one possible $\m \in \Ap(S,E)$.

\begin{lemma}\label{Mij}
Let $S$ be nearly Gorenstein  of embedding dimension $d+3$, and take some indices $s$, $i$, and $j$ such that $1\leq s\leq d$ and $1\leq i<j\leq3$. Then, there exists at most one element $\m\in\max_{\preceq_S}\Ap(S,E)$ such that $T(\a_s)+\a_s=\m+\lambda_i\a_{d+i}+\lambda_j\a_{d+j}$ with $\lambda_i,\lambda_j\in\NN\setminus\{0\}$.
\end{lemma}
\begin{proof}
	Let  $\m,\n \in\max_{\preceq_S}\Ap(S,E)$ such that
	\[
	T(\a_s)+\a_s=\m+\lambda_i\a_{d+i}+\lambda_j\a_{d+j}=\n+\mu_i\a_{d+i}+\mu_j\a_{d+j},
	\]
	with $\lambda_i,\lambda_j,\mu_i,\mu_j\in\NN\setminus\{0\}$. By Lemma~\ref{lll}, there exist $\w_i,\w_j,\w'_i,\w'_j\in\Ap(S,E)$ such that 
	$\m+\a_{d+i}=\w_i+\a_s, \m+\a_{d+j}=\w_j+\a_s, \n+\a_{d+i}=\w'_i+\a_s, \n+\a_{d+j}=\w'_j+\a_s$. 	
	Note that $\w_i-\a_{d+i}, \w'_i-\a_{d+i}, \w_j-\a_{d+j}, \w'_j-\a_{d+j}$ are not in $S$. Therefore,
	$\w_i=h_j\a_{d+j}+h_k\a_{d+k}$, $\w_j=f_i\a_{d+i}+f_k\a_{d+k}$, $\w'_i=h'_j\a_{d+j}+h'_k\a_{d+k}$ and $\w'_j=f'_i\a_{d+i}+f'_k\a_{d+k}$ for some $h_j,h_k,h'_j,h'_k,f_i,f_k,f'_i,f'_k\in\NN$.
	Then 
	\begin{eqnarray*}
		(h_j+1)\a_{d+j}+h_k\a_{d+k}&=&(f_i+1)\a_{d+i}+f_k\a_{d+k}\\
		(h'_j+1)\a_{d+j}+h'_k\a_{d+k}&=&(f'_i+1)\a_{d+i}+f'_k\a_{d+k}.
	\end{eqnarray*}
	Without loss of generality, we assume that $h_k\leq f_k$. Then 
	\[
	(h_j+1)\a_{d+j}=(f_i+1)\a_{d+i}+(f_k-h_k)\a_{d+k}
	\]	 
	implies that $\a_{d+j}$ belongs to the cone generated by $\a_{d+i}$ and $\a_{d+k}$. If $h'_k>f'_k$, then 
	$(h_j'+1)\a_{d+j}+(h'_k-f'_k)\a_{d+k}=(f_i+1)\a_{d+i}$,	 
	will put $\a_{d+i}$ in the cone generated by $\a_{d+j}$ and $\a_{d+k}$, which implies that all three vectors $\a_{d+i}, \a_{d+j}, \a_{d+k}$ are on the same line passing through the origin, a  contradiction by Corollary~\ref{line}. Thus, $h'_k\leq f'_k$. 
	If $h_j=h'_j$, then $\w_i$ and $\w'_i$ are comparable with respect to $\preceq_S$ which means that $\m$ and $\n$ are comparable, a contradiction. So, assume without loss of generality that $h_j<h'_j$. Then 
	\begin{eqnarray*}
		\w'_i&=&(h_j+1)\a_{d+j}+(h'_j-h_j-1)\a_{d+j}+h'_k\a_{d+k}\\&=&(f_i+1)\a_{d+i}+(f_k-h_k+h'_k)\a_{d+k}+(h'_j-h_j-1)\a_{d+j},
	\end{eqnarray*}
	a contradiction since $\w'_i-\a_{d+i}\notin S$.	 
\end{proof}

The previous lemma already implies that the type of $S$ is bounded. Indeed, for each $\m \in \Ap(S,E)\setminus \{T(\a_1)\}$ we know that 
$T(\a_1)+\a_1=\m+\sum^3_{s=1}\lambda_s\a_{d+s}$ with $\lambda_s \in \NN$ for $s=1,2,3$.
Of course, at least one $\lambda_s$ has to be positive. If only one is positive, then by Lemma \ref{lem-T}(2) it does not depend on $\m$, and therefore there is only one possible $\m$ for every index $s=1,2,3$. By the previous lemma, there are at most $6$ possible $\m \in \Ap(S,E)$ with exactly two positive $\lambda_s$. Hence, counting also $T(\a_1)$, this means that the type of $S$ is at most $10$. In order to reduce this bound, in the next lemmas we show that all these elements cannot exist at the same time.

\begin{lemma}\label{maxij}
	Let $S$ be nearly Gorenstein  of embedding dimension $d+3$ and 
	\[
	T(\a_s)+\a_s=\n_1+\mu_i\a_{d+i}=\n_2+\mu_j\a_{d+j},
	\]
	for  some $1\leq s\leq d$, $\n_1,\n_2\in\max_{\preceq_S}\Ap(S,E)$,  $\mu_i,\mu_j\in\NN\setminus\{0\}$ with  $1\leq i<j\leq3$. Then the following statements hold.
	\begin{enumerate}
		\item  $T(\a_s)-\lambda_i\a_{d+i}-\lambda_j\a_{d+j}\in S$ where $\lambda_t=\mu_t-1=\max\{l \ ; \ T(\a_s)-l\a_{d+t}\in S\}$ for $t=i,j$. 
		\item Let	$T(\a_s)+\a_s=\m+\sum^3_{t=1}l_t\a_{d+t}$ where  $\m\in\max_{\preceq_S}\Ap(S,E)$ and $l_t\in\NN$.
		Then $\{t \ ; \ 1\leq t\leq3 , l_t\neq0\}\neq\{i,j\}$.
	\end{enumerate}
\end{lemma}
\begin{proof}
	
	Let $\{i,j,k\}=\{1,2,3\}$ and $s=1$, for simplicity. 
	
	(1) Note that $\mu_t-1=\max\{l \ ; \ T(\a_1)-l\a_{d+t}\in S\}$ for $t=i,j$, by Lemma~\ref{lem-T}(2).  Lemma~\ref{lll} implies that 
	$\n_1+\a_{d+i}=\w+\a_s$ and $\n_2+\a_{d+j}=\w'+\a_s$ for some $\w,\w'\in\Ap(S,E)$. Note that $\w-\a_{d+i}$ and $\w-\a_{d+j}$ are not in $S$. So,
	$\w=h_j\a_{d+j}+h_k\a_{d+k}$ and $\w'=h'_i\a_{d+i}+h_k\a_{d+k}$ for some $h_j,h_k,h'_i.h'_k\in\NN$. 
	Then
	\begin{equation}\label{exp}
	T(\a_s)=(\mu_i-1)\a_{d+i}+h_j\a_{d+j}+h_k\a_{d+k}=(\mu_j-1)\a_{d+j}+h'_i\a_{d+i}+h'_k\a_{d+k}.
	\end{equation}
		Consequently, 
	\[
	(\mu_i-1-h'_i)\a_{d+i}+h_k\a_{d+k}=(\mu_j-1-h_j)\a_{d+j}+h'_k\a_{d+k}.
	\]
	 Note that, $h'_i\leq\mu_i-1$ and $h_j\leq\mu_j-1$. If one of them is an equality, we have done. Assume that it is not the case.  
	If $h_k\neq h'_k$, then either $\mu_i\a_{d+i}-(\a_{d+i}+\a_{d+j}+\a_{d+k})\in S$ or  $\mu_j\a_{d+i}-(\a_{d+i}+\a_{d+j}+\a_{d+k})\in S$, both contradict  Corollary~\ref{j=0}. Therefore, $h_k=h'_k$ and  so $(\mu_i-1-h'_i)\a_{d+i}=(\mu_j-1-h_j)\a_{d+j}$. Then
	\[
	T(\a_s)+\a_s=\n_1+(1+h'_i)\a_{d+i}+(\mu_j-1-h_j)\a_{d+j}=\n_2+(1+h_j)\a_{d+j}+(\mu_i-1-h'_i)\a_{d+i}.
	\]
	Now, Lemma~\ref{Mij} implies that $h_j=\mu_j-1$ or $h'_i=\mu_i-1$, which along with (\ref{exp}) yields the result.
	
	\medskip

	(2) By the statement (1)  
	\[
	T(\a_s)=\lambda_i\a_{d+i}+\lambda_j\a_{d+j}+\mu\a_{d+k},
	\]
	for some $\mu\in\NN$. Let $L=\{t \ ; \ 1\leq t\leq3 \ , \ l_t\neq0\}$.  Assume on the contrary that  $L=\{i,j\}$. Then 
	\begin{equation}\label{ij}
	T(\a_s)+\a_s=\m+l_i\a_{d+i}+l_j\a_{d+j},
	\end{equation}
	with $l_i,l_j\in\NN\setminus\{0\}$. 
By Lemma~\ref{lll},
	$\m+\a_{d+i}=\w+\a_s$ and $\m+\a_{d+j}=\w'+\a_s$ for some $\w,\w'\in\Ap(S,E)$. 	
	Note that $\w-\a_{d+i}$ and $\w'-\a_{d+j}$ are not in $S$. Therefore, $\w=h_j\a_{d+j}+h_k\a_{d+k}$ and $\w'=h_i'\a_{d+i}+h'_k\a_{d+k}$ for some $h_j,h_k,h'_i,h'_k\in\NN$.
	Then  we get
	\begin{eqnarray*}
		T(\a_s)&=&\lambda_i\a_{d+i}+\lambda_j\a_{d+j}+\mu\a_{d+k}\\
		&=&(l_i-1)\a_{d+i}+(h_j+l_j)\a_{d+j}+h_k\a_{d+k}\\
		&=&(h'_i+l_i)\a_{d+i}+(l_j-1)\a_{d+j}+h'_k\a_{d+k}.
	\end{eqnarray*}
	Note that $\lambda_i\geq l_i-1 , h'_i+l_i$ and $\lambda_j\geq l_j-1 , h_j+l_j$. Therefore, $h_k, h'_k\geq\mu$ and 
	\[
	(\lambda_i-l_i+1)\a_{d+i}+(\lambda_j-h_j-l_j)\a_{d+j}=(h_k-\mu)\a_{d+k},
	\]
	\[
	(\lambda_i-h'_i-l_i)\a_{d+i}+(\lambda_j-l_j+1)\a_{d+j}=(h'_k-\mu)\a_{d+k}.
	\]
	Since $\w-\a_{d+i}$ and $\w'-\a_{d+j}$ are not in $S$, we should have $h_k=h'_k=\mu$, which implies $\lambda_i=l_i-1=l_i+h'_i$, a contradiction. 	
\end{proof}

\begin{lemma}\label{M123}
	Let $S$ be  nearly Gorenstein  of embedding dimension $d+3$ such that
	\[
 \w=\m_1+\mu_1\a_{d+1}=\m_2+\mu_2\a_{d+2}=\m_3+\mu_3\a_{d+3},
	\]
	where 	 $\m_1,\m_2,\m_3$  are three different elements in  $\max_{\preceq_S}\Ap(S,E)$ and  $\mu_1,\mu_2,\mu_3\in\NN$. Then  $\w\neq T(\a_i)+\a_i$ for  any $1\leq i\leq d$.
\end{lemma}
\begin{proof}
	Assume on the contrary that $\w=T(\a_i)+\a_i$ for some $1\leq i\leq d$. By Lemma~\ref{maxij}(1), 
\begin{eqnarray*}
T(\a_i)&=&\lambda_1\a_{d+1}+\lambda_2\a_{d+2}+\lambda'_3\a_{d+3}\\
&=&\lambda_1\a_{d+1}+\lambda'_2\a_{d+2}+\lambda_3\a_{d+3}	\\
&=&\lambda'_1\a_{d+1}+\lambda_2\a_{d+2}+\lambda_3\a_{d+3},
\end{eqnarray*}
where $\lambda'_1,\lambda'_2,\lambda'_3\in\NN$ and $\lambda_t=\mu_t-1=\max\{l \ ; \ T(\a_i)-l\a_{d+t}\in S\}$ for $t=i,j$. If $\lambda_t>\lambda'_t$ for some $1\leq t\leq 3$, then 
\[
(\lambda_3-\lambda'_3)\a_{d+3}=(\lambda_2-\lambda'_2)\a_{d+2}=(\lambda_1-\lambda'_1)\a_{d+1},
\]
a contradiction with Corollary~\ref{line}. Therefore, $\lambda_t=\lambda'_t$ for $t=1,2,3$ which means that $T(\a_i)$ has a unique expression. By Proposition~\ref{uni-ex-r}, we get  $\typ(S)=4$ which implies that  $d\geq2$, by  
\cite[Theorem~2.4]{Moscariello-Strazzanti}.
Thus,  $\max_{\preceq_S}\Ap(S,E)=\{T(\a_i), \m_1,\m_2,\m_3\}$. Assume, without loss of generality, that $T(\a_j)=\m_1$ for some $1\leq j\neq i\leq d$. Then 
\begin{equation}\label{Tai}
T(\a_i)+\a_i=T(\a_j)+\mu_1\a_{d+1}=\m_2+\mu_2\a_{d+2}=\m_3+\mu_3\a_{d+3},
\end{equation}
\[
T(\a_j)+\a_j=\m_2+l_1\a_{d+1}+l_3\a_{d+3}=\m_3+h_1\a_{d+1}+h_2\a_{d+2},
\]
where $l_1,l_3,h_1,h_2\in\NN$, by Lemma~\ref{lem-T}(1). Note that
\begin{equation}\label{lh}
h_1\a_{d+1}+(\mu_2+h_2)\a_{d+2}=l_1\a_{d+1}+(\mu_3+l_3)\a_{d+3},
\end{equation}
 from Lemma~\ref{lem-T}(3). By this non-trivial relation it follows that $\a_{d+1},\a_{d+2},\a_{d+3}$ belong to a two dimensional cone, which implies $d=2$  by Proposition~\ref{hplane}. Without loss of generality, let $h_1\leq l_1$. Then $\a_{d+2}$ belongs to the cone generated by $\a_{d+1}$ and $\a_{d+3}$. 
We may order the generating vectors by their slopes as the following 
\[
\slop(\a_1)\leq\slop(\a_{d+1})\leq\slop(\a_{d+2})\leq\slop(\a_{d+3})\leq\slop(\a_2),
\]
where $\{1,2\}=\{i,j\}$. 
By Lemma~\ref{lem-T}(3),
\begin{equation*}\label{j}
(\mu_1+h_1)\a_{d+1}+h_2\a_{d+2}=\mu_3\a_{d+3}+\a_j,
\end{equation*}
Which 
 implies $j=1$ and so $i=2$.  If  $\a_{d+1}$ and $\a_{d+2}$ have the same slope, then all three vectors are on the same line passing through the origin of coordinates  by (\ref{lh}), which makes a   contradiction by  Corollary~\ref{line}. Thus,
 \begin{equation}\label{slopp}
 \slop(\a_1)\leq\slop(\a_{d+1})<\slop(\a_{d+2})\leq\slop(\a_{d+3})\leq\slop(\a_2).
 \end{equation}
 If $T(\a_j)$ has two different expressions, then there exists $f_2\in\NN\setminus\{0\}$ such that $f_2\a_{d+2}\preceq_S T(\a_j)$ and
 $f_2\a_{d+2}=f_1\a_{d+1}+f_3\a_{d+3}$ for some $f_1,f_3\in\NN$. 
Then $f_2\geq\mu_2$, which along with (\ref{Tai}) implies that $\m_2-\mu_1\a_{d+1}\in S$. As $T(\a_i)=\sum^3_{t=1}(\mu_t-1)\a_{d+t}$, looking again at (\ref{Tai}), we get 
\[
(\mu_3-1)\a_{d+3}+\a_2=\m_2-(\mu_1-1)\a_{d+1}+\a_{d+2}.
\]
 As $\a_{d+1}$ appears in an expression of the right hand side of the above equation, it makes a contradiction by the order of slopes in (\ref{slopp}). Therefore, $T(\a_j)$ has a unique expression and so  $0\in\{l_1,l_3\}\cap\{h_1,h_2\}$, by Lemma~\ref{uni-ex-li}. If $l_1=0$, then $h_1=0$ and the equation (\ref{lh}) implies that $\a_{d+2}$ and $\a_{d+3}$ have the same slope, a contradiction.
Thus, $l_3=0$. As $\m_2$ and $\m_3$ are not comparable with respect to $\preceq_S$, we have $h_2\neq 0$ and so $h_1=0$. Therefore,
\[
T(\a_1)+\a_1=T(\a_2)+e_3\a_{d+3}=\m_2+l_1\a_{d+1}=\m_3+h_2\a_{d+2},
\]
for some $e_3\in\NN$ and
\[
T(\a_2)+\a_2=T(\a_1)+\mu_1\a_{d+1}=\m_2+\mu_2\a_{d+2}=\m_3+\mu_3\a_{d+3}.
\]
Note that $T(\a_1)=(l_1-1)\a_{d+1}+(h_2-1)\a_{d+2}+(e_3-1)\a_{d+3}$ by Lemma~\ref{lem-T}(2), and  $T(\a_2)+e_3\a_{d+3}+\mu_2\a_{d+2}=T(\a_1)+\mu_1\a_{d+1}+l_1\a_{d+1}$. Thus,
\[
(\mu_1-1)\a_{d+1}+(2\mu_2-1)\a_{d+2}+(\mu_3-1+e_3)\a_{d+3}
=(\mu_1+2l_1-1)\a_{d+1}+(h_2-1)\a_{d+2}+(e_3-1)\a_{d+3},
\]
which implies 
\[
(2\mu_2-1)\a_{d+2}+\mu_3\a_{d+3}=2l_1\a_{d+1}+(h_2-1)\a_{d+2}.
\]
This equation is impossible by (\ref{slopp}).  
\end{proof}

The previous lemma allows us to prove that $\typ(S) \leq 3$ when some $T(\a_i)$ has a unique expression. We will make use of this fact later.

\begin{corollary}\label{uni-ex}
	If $S$ is nearly Gorenstein of embedding dimension $d+3$ and $T(\a_i)$ has a  unique expression for some $1\leq i\leq d$, then $\typ(S)\leq 3$. 
\end{corollary}
\begin{proof}
Note that $\typ(S)\leq 4$, by 	Proposition~\ref{uni-ex-r}. Assume on the contrary that $\typ(S)=4$ and let $\max_{\preceq_S}\Ap(S,E)=\{T(\a_1),\m_1,\m_2,\m_3\}$. Then by Lemma~\ref{uni-ex-li},
\[
T(\a_1)+\a_1=\m_1+l_i\a_{d+i}=\m_2+l_j\a_{d+j}=\m_3+l_k\a_{d+k},
\]
where $\{i,j,k\}=\{1,2,3\}$ and $l_i,l_j,l_k\in\NN$. This is a contradiction by Lemma~\ref{M123}. 
\end{proof}

In order to prove that the type of $S$ is always at most three, we need other two lemmas. 

\begin{lemma}\label{Mijk}
		Let $S$ be nearly Gorenstein  of embedding dimension $d+3$ with $d\geq2$ and $1\leq s\leq d$. Let $\m_1,\m_2\in\max_{\preceq_S}\Ap(S,E)$ such that 
		\[
		T(\a_s)+\a_s=\m_1+f_i\a_{d+i}+f_j\a_{d+j}=\m_2+g_i\a_{d+i}+g_k\a_{d+k},
		\]
		where $\{i,j,k\}=\{1,2,3\}$ and $f_i,f_j,g_i,g_k\in\NN\setminus\{0\}$. 
		Then:
		\begin{enumerate}
			\item $T(\a_s)-\a_{d+t}\in S$, for $t=1,2,3$;
			\item $\a_{d+i}$ is an  interior point  of the cone generated by $\a_{d+j}$ and $\a_{d+k}$;
				\item $T(\a_s)+\a_s\neq\n+l\a_{d+i}$, for any $\n\in\max_{\preceq_S}\Ap(S,E)$ and $l\in\NN$;
			\item 	$T(\a_t)\in\{\m_1,\m_2\}$ for $t\in\{1,\dots,d\}\setminus\{s\}$.
		\end{enumerate}
\end{lemma}
\begin{proof}
 By Lemma~\ref{lll},  
$\m_1+\a_{d+i}=h_j\a_{d+j}+h_k\a_{d+k}+\a_s$ and $\m_1+\a_{d+j}=h_i\a_{d+i}+h'_k\a_{d+k}+\a_s$ for some $h_j,h_i,h_k,h'_k\in\NN$.  
Then 
\[
T(\a_s)=(f_i-1)\a_{d+i}+(h_j+f_j)\a_{d+j}+h_k\a_{d+k}=(h_i+f_{i})\a_{d+i}+(f_j-1)\a_{d+j}+h'_k\a_{d+k},
\]
and	consequently
\begin{equation}\label{hh'1}
(h_j+1)\a_{d+j}+h_k\a_{d+k}=(h_i+1)\a_{d+i}+h'_k\a_{d+k}.
\end{equation}
By a similar argument, we get
$\m_2+\a_{d+i}=e_j\a_{d+j}+e_k\a_{d+k}+\a_s$ and $\m_2+\a_{d+k}=e_i'\a_{d+i}+e'_j\a_{d+j}+\a_s$ for some  $e_j,e_k,e'_i,e'_j\in\NN$.  Then
\[
T(\a_s)=(g_{i}-1)\a_{d+i}+(e_k+g_{k})\a_{d+k}+e_j\a_{d+j}=
(e'_i+g_{i})\a_{d+i}+(g_k-1)\a_{d+k}+e'_j\a_{d+j},
\]
implies the statement (1), 
and	
\begin{equation}\label{ee'1}
(e_k+1)\a_{d+k}+e_j\a_{d+j}=(e_i'+1)\a_{d+i}+e'_j\a_{d+j}.
\end{equation}	
	If $h_k\leq h'_k$, then by (\ref{hh'1}) it follows that $\a_{d+j}$ belongs to the cone generated by $\a_{d+i}$ and $\a_{d+k}$. By (\ref{ee'1}), we have either $\a_{d+k}$ belongs to the cone of $\a_{d+i}$ and $\a_{d+j}$ or 
$\a_{d+i}$ belongs to the cone generated by  $\a_{d+k}$ and $\a_{d+j}$, both cases will imply that $\a_{d+i},\a_{d+j},\a_{d+k}$ are on the same line passing through the origin, which contradicts Corollary~\ref{line}.
So that $h_k>h'_k$. Then  $\a_{d+i}$ belongs to the interior of the cone generated by  $\a_{d+k}$ and $\a_{d+j}$ from (\ref{hh'1}), which is the subject of (2).

Note that $h_i+1\leq h_i+f_i\leq\lambda_i$, where $\lambda_i=\max\{l\ ; \ T(\a_s)-l\a_{d+i}\in S\
\}$. If  $T(\a_s)+\a_s=\n+l\a_{d+i}$, for some $\n\in\max_{\preceq_S}\Ap(S,E)$ and $l\in\NN$, then  by Lemma~\ref{lem-T}(2) and the  equation (\ref{hh'1}),
\begin{equation*}
T(\a_s)+\a_s=\m_i+(\lambda_i+1)\a_{d+i}=\m_i+(\lambda_i-h_i)\a_{d+i}+(h_j+1)\a_{d+j}+(h_k-h'_k)\a_{d+k},
\end{equation*}
which is a contradiction by Corollary~\ref{j=0}. Thus the statement (3), is obtained. 

In order to prove (4), assume on the contrary that $T(\a_t)\notin\{\m_1, \m_2\}$, then 
\[T(\a_t)+\a_t=\m_2+\mu_j\a_{d+j}=\m_1+\mu_k\a_{d+k},\]
for some $\mu_j,\mu_k\in\NN$, by Lemma~\ref{lem-T}(1). Now, using Lemma~\ref{lem-T}(3), we get 
\[g_i\a_{d+i}+g_k\a_{d+k}+\mu_k\a_{d+k}=f_i\a_{d+i}+f_j\a_{d+j}+\mu_j\a_{d+j}.\]
Consequently, 
\[ (g_i-f_i)\a_{d+i}+(g_k+\mu_k)\a_{d+k}=(f_j+\mu_j)\a_{d+j},
\]
which along with the statement (2),  implies that  $\a_{d+i}, \a_{d+j},\a_{d+k}$ are on the same line passing through the origin, a contradiction by Corollary~\ref{line}. 
\end{proof}

\begin{lemma}\label{lem-MM}
	Let $S$ be nearly Gorenstein  of embedding dimension $d+3$ with $d\geq2$  and let 	$$M^1_i=\{\m\in\max_{\preceq_S}\Ap(S,E) \ ; \ T(\a_1)+\a_1=\m+\a_{d+i}  \},$$
		$$M^2_i=\{\m\in\max_{\preceq_S}\Ap(S,E) \ ; \ T(\a_1)+\a_1=\m+\lambda\a_{d+i} \text{ with } \lambda\geq2 \},$$
	$$M_{i,j}=\{\m\in\max_{\preceq_S}\Ap(S,E) \ ; \ T(\a_1)+\a_1=\m+\lambda\a_{d+i}+\mu\a_{d+j} \text{ with } \lambda,\mu\in\NN\setminus\{0\} \},$$
	for $1\leq i,j\leq 3$ and $i \neq j$. Then the following statements hold, where $\{i,j,k\}=\{1,2,3\}$.
	\begin{enumerate}
		\item Each of $M^1_i, M^2_i$ and $M_{i,j}$ has at most one element.
			\item If $M^1_i$ is not empty, then $M^2_i=M_{i,j}=M_{i,k}=\emptyset$.
				\item  If $M_{i,j}$  and $M_{i,k}$ are not  empty, then  $M_{j,k}=M^1_t=M^2_t=\emptyset$ for $t=1,2,3$. 	
				\item If $M^2_i$ and $M^2_j$ are  two different non-empty sets, then  $M^2_k\cup M^1_k\subset M^2_i\cup M^2_j$ and
				$M^1_t=M_{i,j}=M_{i,k}=M_{j,k}=\emptyset$, for $t=i,j$.

		
	\end{enumerate}
\end{lemma}
\begin{proof}
	(1) If $\m,\n\in M^1_i$, then $T(\a_1)+\a_1=\m+\lambda\a_{d+i}=\n+\mu\a_{d+i}$ for some $\lambda,\mu\in\NN$, which implies $\m=\n$, since $\m$ and $\n$ are maximal elements. Now, Lemma~\ref{Mij} completes the proof.

\medskip

Assume that $M^1_t\subseteq\{\n_t\}$, $M^2_t\subseteq\{\m_t\}$ 
 and $M_{t,s}\subseteq\{\m_{t,s}\}$ for $1\leq t,s\leq3$. In the case that $M^2_t$ is not empty, we let 
	\begin{equation}\label{eq-i1}
	T(\a_1)+\a_1=\m_t+\mu_t\a_{d+t},
	\end{equation}
	for some integer $\mu_t\geq2$ and in the case that $M_{t,s}$ is not empty,  let 
	\begin{equation}\label{eq-ij}
	T(\a_1)+\a_1=\m_{t,s}+f_{t,s}\a_{d+t}+g_{t,s}\a_{d+s},
	\end{equation}
	for some $f_{t,s}, g_{t,s}\in\NN\setminus\{0\}$.

	\medskip

(2) As $T(\a_1)+\a_1=\n_i+\a_{d+i}$, if any of $M^2_i, M_{i,j}, M_{i,k}$ is non-empty, then  $\n_i\in\{\m_i+(\mu_i-1)\a_{d+i}, \m_{i,j}+(f_{i,j}-1)\a_{d+i}+g_{i,j}\a_{d+j},\m_{i,k}+(f_{i,k}-1)\a_{d+i}+g_{i,k}\a_{d+k}\}$, a contradiction.

	\medskip
	
(3)  If $M_{j,k}\neq\emptyset$, then by Lemma~\ref{Mijk}(2),	each of the vectors $\a_{d+i},\a_{d+j},\a_{d+k}$ belongs to the cone generated by two others, which means that they  are on the same line passing through the origin, a contradiction by Corollary~\ref{line}. Thus, $M_{j,k}=\emptyset$. As $T(\a_1)-\a_{d+t}\in S$, for $t=i,j,k$, by Lemma~\ref{Mijk}(1), we get $M^1_i=M^1_j=M^1_k=\emptyset$ along with Lemma~\ref{lem-T}(2). We also have $M^2_i=\emptyset$, by Lemma~\ref{Mijk}(3).  
		Note that 
	\begin{equation}\label{hi}
	h_i\a_{d+i}=h_j\a_{d+j}+h_k\a_{d+k},
	\end{equation}
	for some $h_i,h_j,h_k\in\NN\setminus\{0\}$, 
 by Lemma~\ref{Mijk}(2). It  follows $d=2$ by Proposition~\ref{hplane}. 
	We may order the generating vectors by their slopes as the following 
	\begin{equation}\label{slop}
	\slop(\a_p)\leq\slop(\a_{d+k})\leq\slop(\a_{d+i})\leq\slop(\a_{d+j})\leq\slop(\a_q),
	\end{equation}
	where $\{p,q\}=\{1,2\}$. 
	Note that $T(\a_2)\in\{\m_{i,k}, \m_{i,j}\}$ by Lemma~\ref{Mijk}(4). Without loss of generality, assume that $T(\a_2)=\m_{i,j}$. By Lemma~\ref{lem-T}(1)
\begin{eqnarray*}
	T(\a_1)+\a_1&=&T(\a_2)+f_{i,j}\a_{d+i}+g_{i,j}\a_{d+j}=\m_{i,k}+f_{i,k}\a_{d+i}+g_{i,k}\a_{d+k},\\ 
	T(\a_2)+\a_2&=&\m_{i,k}+f_j\a_{d+j},
\end{eqnarray*}
where  $f_j\in\NN$, and by Lemma~\ref{lem-T}(5) we get
	\[
	f_{i,j}\a_{d+i}+(g_{i,j}+f_j)\a_{d+j}=f_{i,k}\a_{d+i}+g_{i,k}\a_{d+k}+\a_2,
	\]
	which implies $q=2$ and so $p=1$. If  $M^2_j\cup M^2_k\neq\emptyset$, then either $T(\a_1)+\a_1=\m_j+\mu_j\a_{d+j}$ or $T(\a_1)+\a_1=\m_k+\mu_k\a_{d+k}$, for some $\mu_j,\mu_k\in\NN$. Comparing with the above equations of $T(\a_1)+\a_1$, and the fact that $\a_{d+i},\a_{d+j},\a_{d+k}$  are not on the same line passing through the origin, we get that $T(\a_1), T(\a_2), \m_{i,k},\m_t$, with $t\in\{j,k\}$, are four different elements. Therefore, $T(\a_2)$ does not have a unique expression, by Corollary~\ref{uni-ex}. Let $h_i$ be the minimum positive integer satisfied the equation~\ref{hi}. Then $T(\a_2)$ has an expression
	$T(\a_2)=\sum^3_{t=1}l_t\a_{d+t}$ with $l_t\in\NN$ and $l_i\geq h_i$. As 
	\[
		T(\a_1)+\a_1=T(\a_2)+f_{i,j}\a_{d+i}+g_{i,j}\a_{d+j},
	\]
	there exists $e_i,e_k\in\NN$ such that $T(\a_2)+\a_{d+j}=e_i\a_{d+i}+e_k\a_{d+k}+\a_1$, by Lemma~\ref{lll}. Then $e_i<h_i\leq l_i$, as $T(\a_2)-\a_1\notin S$. Therefore,
	\[
	(l_i-e_i)\a_{d+i}+(l_j+1)\a_{d+j}+l_k\a_{d+k}=e_k\a_{d+k}+\a_1.
	\]
Since $l_i-e_i>0$ and $l_j+1>0$, we have $e_k-l_k>0$. 
Considering the order of slopes in (\ref{slop}), the above equation implies that $\a_{d+i},\a_{d+j},\a_{d+k}$ have the same slope, a contradiction with Corollary~\ref{line}. Therefore, $M^2_j\cup M^2_k=\emptyset$.

\medskip

(4)  Note that $M^2_k\cup M^1_k\subset M^2_i\cup M^2_j$  by Lemma~\ref{M123}  and  $M_{i,j}=\emptyset$ by Lemma~\ref{maxij}(2). 	 Since $M^2_i$ is not empty, $T(\a_1)-\a_{d+i}\in S$ and so $M^1_i=\emptyset$, by Lemma~\ref{lem-T}(1).
		Now, assume on the contrary that $M_{i,k}$ is not empty.  Then 
		\begin{equation}\label{ik}
		T(\a_1)+\a_1=\m_{i,k}+f_{i,k}\a_{d+i}+g_{i,k}\a_{d+k}.
		\end{equation}
		By Lemma~\ref{lll},   
		$\m_{i,k}+\a_{d+i}=\w+\a_1$ and $\m_{i,k}+\a_{d+k}=\w'+\a_1$ for some $\w,\w'\in\Ap(S,E)$. 	
		Note that $\w-\a_{d+i}$ and $\w'-\a_{d+k}$ are not in $S$. Therefore, $\w=h_j\a_{d+j}+h_k\a_{d+k}$ and $\w'=h_i'\a_{d+i}+h'_j\a_{d+j}$ for some $h_j,h_k,h'_i,h'_j\in\NN$.  Then using also Lemma~\ref{maxij}(2) we get
		\begin{eqnarray*}
			T(\a_1)&=&\lambda_i\a_{d+i}+\lambda_j\a_{d+j}+\mu\a_{d+k}\\
			&=&(f_{i,k}-1)\a_{d+i}+h_j\a_{d+j}+(g_{i,k}+h_k)\a_{d+k}\\
			&=&(h'_i+f_{i,k})\a_{d+i}+h'_j\a_{d+j}+(g_{i,k}-1)\a_{d+k},
		\end{eqnarray*}
		where $\lambda_t=\max\{l \ ; \ T(\a_1)-l\a_{d+t}\in S\}$ for $t=i,j$.
	Then 
				$\lambda_i\geq f_{i,k}+h'_i>f_{i,k}-1$ and $\lambda_j\geq h_j , h'_j$. Thus, $\mu\leq g_{i,k}-1<g_{i,k}+h_k$ and 
	the equation
		\begin{equation}\label{a}
		(g_{i,k}+h_k-\mu)\a_{d+k}=(\lambda_i-f_{i,k}+1)\a_{d+i}+(\lambda_j-h_j)\a_{d+j},
		\end{equation}
	shows 		that $\a_{d+1},\a_{d+2},\a_{d+3}$ belong to a two dimensional cone, which implies $d=2$  by Proposition~\ref{hplane}. 
		We may order the generating vectors by their slopes as the following 
	\begin{equation}\label{slp}
		\slop(\a_s)\leq\slop(\a_{d+i})\leq\slop(\a_{d+k})\leq\slop(\a_{d+j})\leq\slop(\a_t),
	\end{equation}
		where $\{s,t\}=\{1,2\}$. 
	Note that $\m_{i,k},\n_i,\n_j$ are three different elements. If $T(\a_2)\notin\{\n_i,\n_j\}$, then 
		\[
		T(\a_2)+\a_2=\n_i+f_j\a_{d+j}+f_k\a_{d+k}=\n_j+g_i\a_{d+i}+g_k\a_{d+k},
		\]
		by Lemma~\ref{lem-T}(1), and 
		\[
		(\mu_i+g_i)\a_{d+i}+g_k\a_{d+k}=(\mu_j+f_j)\a_{d+j}+f_k\a_{d+k},
		\]
		by Lemma~\ref{lem-T}(3). But the recent equation, along with (\ref{a}), implies that $\a_{d+i},\a_{d+j},\a_{d+k}$ are on the same line passing through the origin of coordinates, contradiction by Corollary~\ref{line}. Therefore, $T(\a_2)\in\{\n_i,\n_j\}$ and 
		\[
		T(\a_2)+\a_2=\m_{i,k}+(\lambda_j+1)\a_{d+j},
		\]
		by Lemma~\ref{lem-T}(1).  
				If $T(\a_2)=\n_i$, then 
		\[
		T(\a_2)+\a_2=\n_j+g_i\a_{d+i}+g_k\a_{d+k},
		\]	
		for some $g_i,g_k\in\NN$, by Lemma~\ref{lem-T}(1), and
		\[
		(\mu_i+g_i)\a_{d+i}+(\mu_k+g_k)\a_{d+k}=\mu_j\a_{d+j}+(\lambda_j-1)\a_{d+j}.
		\]
		by Lemma~\ref{lem-T}(3). But the recent equation, along with (\ref{a}), implies that $\a_{d+i},\a_{d+j},\a_{d+k}$ have the same slope, contradiction by Corollary~\ref{line}. Therefore, $T(\a_2)\neq\n_i$ and	so $T(\a_2)=\n_j$ and 
		\begin{equation}\label{Ta1}
		T(\a_1)+\a_1=T(\a_2)+(\lambda_j+1)\a_{d+j}=\n_1+(\lambda_i+1)\a_{d+i}=\m_{i,k}+f_{i,k}\a_{d+i}+g_{i,k}\a_{d+k},
		\end{equation}
		\begin{equation}\label{Ta2}
		T(\a_2)+\a_2=T(\a_1)+\sum^3_{t=1}q_t\a_{d+t}=\n_i+f_j\a_{d+j}+f_k\a_{d+k}=\m_{i,k}+\mu_j\a_{d+j},
		\end{equation}
		where $q_1,q_2,q_3,\mu_j\in\NN$. Note that $f_k\neq0$, because otherwise $\m_{i,k}$ and $\n_i$ will be comparable with respect $\preceq_S$. 
If $q_i=0$, then we have the following cases:
\begin{enumerate}
	\item If $f_j>0$, then either $q_j=0$ or $q_k=0$, by Lemma~\ref{Mij}. As $T(\a_1)$ is not comparable with $\m_{i,j}$, $q_k>0$. Thus, $q_j=0$ which is contradicts  Lemma~\ref{maxij}(2). 
	\item If $f_j=0$, then $q_j>0$ and $q_k>0$, because $T(\a_1),\n_i,\m_{i,k}$ are not comparable. This contradicts  Lemma~\ref{maxij}(2).
\end{enumerate}		
Therefore, $q_i>0$. By  Lemma~\ref{lem-T}(4), $\sum^3_{t=1}q_t\a_{d+t}+(\lambda_j+1)\a_{d+j}=\a_1+\a_2$. Let $\c=\sum^3_{t=1}q_t\a_{d+t}+(\lambda_j+1)\a_{d+j}-\a_{d+i}$. Then $$\c=(q_i-1)\a_{d+i}+(q_j+\lambda_j+1)\a_{d+j}+q_k\a_{d+k}\in\Ap(S,E),$$ and so $\c\preceq_S\n\in\max_{\preceq_S}\Ap(S,E)$. As  $M^2_k\cup M^1_k\subset M^2_i\cup M^2_j$ and  $M_{i,j}=\emptyset$, 
 $$\max_{\preceq_S}\Ap(S,E)=\{T(\a_1), T(\a_2)=\n_j, \n_i,\m_{i,k}\}\cup M_{j,k}.$$ If $M_{j,k}$ is not empty, then $T(\a_2)+\a_2=\m_{j,k}+p_i\a_{d+i}$ for some $p_i\in\NN$. 
 Note that  
 $\a_1+\a_2=\c+\a_{d+i}\npreceq_S T(\a_1)+\a_1$ and $\a_1+\a_2=\c+\a_{d+i}\npreceq_S T(\a_2)+\a_2$. Therefore, $\c\npreceq_S\m$ for $\m\in\{\n_i, \m_{i,k} , \m_{j,k}\}$. Since $\lambda_j\a_{d+j}\preceq_S \c$, we have 
 $\c\npreceq_S T(\a_1)$. Thus, $\c\preceq_S T(\a_2)=\n_j$. Let $T(\a_2)=\c+\v$.  
 Then $\v-\a_{d+i}\notin S$ and
 \[
 \sum^3_{t=1}q_t\a_{d+t}+(\lambda_j+1)\a_{d+j}-\a_{d+i}+\v+\a_2=T(\a_2)+\a_2=T(\a_1)+\sum^3_{t=1}q_t\a_{d+t},
 \]
 from (\ref{Ta2}). Recall that $T(\a_1)=\lambda_i\a_{d+i}+\lambda_j\a_{d+j}+\mu\a_{d+k}$ and $\v-\a_{d+i}\notin S$. Let $\v=r_j\a_{d+j}+r_k\a_{d+k}$. Then 
 \begin{equation}\label{eq}
(r_j+1)\a_{d+j}+r_k\a_{d+k}+\a_2=(\lambda_i+1)\a_{d+i}+\mu\a_{d+k}.
\end{equation}
By Lemma~\ref{lem-T}(3), we have
\[
(\lambda_j+1+\mu_j)\a_{d+j}=l_i\a_{d+i}+l_k\a_{d+k}+\a_2.
\]
which implies that in (\ref{slp}), $t=j$ and so $s=i$. Now,   equation (\ref{eq}) is in contradiction with (\ref{slp}).   Therefore, $M_{i,k}=\emptyset$. The same argument, replacing $i$ with $j$,  shows that $M_{j,k}=\emptyset$. 
\end{proof}

We are now ready to prove our main theorem.

\begin{theorem}\label{thm}
	Let $S$ be nearly Gorenstein of embedding dimension $d+3$.  If $S$ is not Gorenstein, then   $d\leq\typ(S)\leq3$.
\end{theorem}
\begin{proof}
	By Corollary~\ref{dtyp}, it is enough to show that $\typ(S)\leq3$.
We may assume that $d>1$ by \cite[Theorem~2.4]{Moscariello-Strazzanti} and $\typ(S)\geq 2$. 
Let $M^1_i, M^2_i, M_{i,j}$ be as defined in Lemma~\ref{lem-MM} for $1\leq i,j\leq 3$ and $i\neq j$.  	Then $\max_{\preceq_S}\Ap(S,E)=\{T(\a_1)\}
	\cup(\cup^3_{i=1} M^1_i)\cup(\cup^3_{i=1} M^2_i)\cup
 M_{1,2}\cup M_{1,3}\cup M_{2,3}$, by Corollary~\ref{j=0}.

	Note that each $M^1_i$, $M^2_i$ and $M_{i,j}$ has at most one element by Lemma~\ref{lem-MM}(1). We distinguish the following cases:

	Case 1, $M_{i,j}=M_{i,k}=M_{j,k}=\emptyset$.  If $M^2_i, M^2_j, M^2_k$ are not empty, then  $M^2_i\cup M^2_j\cup M^2_k$ has at most two elements,  by Lemma~\ref{M123}. We may assume that   $M^2_k\subseteq M^2_i\cup M^2_j$. If $M^2_i$ and $M^2_j$ are both non-empty, then by Lemma~\ref{maxij}(1),
	\[
	T(\a_1)=\lambda_i\a_{d+i}+\lambda_j\a_{d+j}+\mu\a_{d+k},
	\]
where $\lambda_t=\max\{l \ ; \ T(\a_1)-l\a_{d+t}\in S\}$ for $t=i,j$. Now, if $M^1_k\neq\emptyset$, then $T(\a_1)-\a_{d+k}\notin S$, by Lemma~\ref{lem-T}(2) and so $\mu=0$. This means that $T(\a_1)$ has a unique expression which implies  $\typ(S)\leq3$ by Corollary~\ref{uni-ex}. So, we may also assume that  $M^1_k=\emptyset$. 
As $M^2_i$ and $M^2_j$ are not empty, we have $M^1_i=M^1_j=\emptyset$, by Lemma~\ref{lem-MM}(2).
If  $M^2_j$ is also empty, then the only possible non-empty sets are $M^2_i=M^2_k$,  $M^1_j$ and  $M^1_k$,   which implies $\typ(S)\leq3$.

Case 2, $M_{i,j}$ and $M_{i,k}$ are not empty. Then $M_{j,k}=M^1_t =M^2_t=\emptyset$ for $t=1,2,3$, by Lemma~\ref{lem-MM}(3).

 Case 3, $M_{i,j}$ is not empty, but $M_{i,k}=M_{j,k}=\emptyset$. 
	 Then $M^1_i=M^1_j=\emptyset$, 
	and $M^2_i\cup M^2_k\cup M^2_j$ has at most one element by Lemma~\ref{lem-MM}(4). 
	Note that $M^2_k$ and $M^1_k$ can not be non-empty at the same time, by Lemma~\ref{lem-MM}(2).
	Therefore, it is enough to show that if $M^2_t$ is not empty for some $t\in\{i,j\}$, then $M^1_k=\emptyset$. Without loss of generality, we assume that $t=i$. Assume on the contrary that
	\begin{equation}\label{1}
	T(\a_1)+\a_1=\m_{i,j}+f_{i,j}\a_{d+i}+g_{i,j}\a_{d+j}=\m_i+\mu_i\a_{d+i}=\n_k+\a_{d+k}.
	\end{equation}
 	By Lemma~\ref{lll},   
	$\m_{i,j}+\a_{d+i}=h_j\a_{d+j}+h_k\a_{d+k}+\a_1$ and $\m_{i,j}+\a_{d+j}=h_i\a_{d+i}+h'_k\a_{d+k}+\a_1$ for some $h_j,h_k,h'_k\in\NN$. 	Then
	\[
	T(\a_1)=(f_{i,j}-1)\a_{d+i}+(h_j+g_{i,j})\a_{d+j}+h_k\a_{d+k}=(h_i+f_{i,j})\a_{d+i}+(g_{i,j}-1)\a_{d+j}+h'_k\a_{d+k}  
	\]
	If $M^1_k\neq\emptyset$, then 	$T(\a_1)-\a_{d+k}\notin S$, by Lemma~\ref{lem-T}(2) and so $h_k=h'_k=0$. In particular,
	\[
	(h_i+1)\a_{d+i}=(h_j+1)\a_{d+j}.
	\]

For $\v\in\{\m_{i,j},\m_i,\n_k\}$, let $[\v]_k=l\a_{d+k}$, where $l=\max\{l \ ; \ \v-l\a_{d+k}\in S\}$. Then, along with (\ref{1}), $$[\m_{i,j}]_k=[\m_i]_k=[\n_k]_k+1,$$ as $\a_{d+i},\a_{d+j},\a_{d+k}$ are not on the same line, by Corollary~\ref{line}.
If $T(\a_2)\notin\{\n_k,\m_{i,j}\}$, then 
\[
T(\a_2)+\a_2=\n_k+g_i\a_{d+i}+g_j\a_{d+j}=\m_{i,j}+g_k\a_{d+k}.
\]
Since $[\m_{i,j}]_k+g_k>[\n_k]_k$, it provides a non-trivial relation between $\a_{d+k}$ and the points on the line passing through $\a_{d+i}$ and $\a_{d+j}$, which means all three points are on the same line,  a contradiction. 

If $T(\a_2)=\n_k$, then 
\[
T(\a_2)+\a_2=\m_{i,j}+g_k\a_{d+k}=\m_i+e_j\a_{d+j}+e_k\a_{d+k}.
\]
Note that $g_k>e_k$. As $[\m_i]_k=[\m_{i,j}]_k$, the above equation provides a non-trivial relation between $\a_{d+k}$ and the points on the line passing through $\a_{d+i}$ and $\a_{d+j}$, which means all three points are on the same line. 

If $T(\a_2)=\m_{i,j}$, then 
\[
T(\a_2)+\a_2=\n_k+g_i\a_{d+i}+g_j\a_{d+j}=\m_i+e_j\a_{d+j}+e_k\a_{d+k}.
\]
As  $\a_{d+k}$ is not on the line passing through $\a_{d+i}$ and $\a_{d+j}$, we get $[\n_k]_k=[\m_i]_k+e_k=[\n_k]_k+1+e_k$, a contradiction. 
\end{proof}

The bounds obtained in the previous theorem are sharp. More precisely, for each possible value $t$ between these bounds there exist nearly Gorenstein semigroups with embedding dimension $d+3$ having type $t$. This is well known for $d=1$, see for instance \cite[Example 2.7]{Moscariello-Strazzanti}. When $d=2$, in Example \ref{d=t=2} we have seen such a semigroup having type $2$. For the two missing cases we provide examples below.

\begin{example} \label{ex1}
Let  $\a_1=(5,0), \a_2=(0,3), \a_3=(3,1), \a_4=(1,2), \a_5=(2,2)$ and let $S$ be the affine semigroup generated by them. The extremal rays of $S$ are $\a_1$ and $\a_2$, and $\max_{\preceq_S}\Ap(S,E)=\{\m_1=(5,10), \m_2=(8,8), \m_3=(7,9)\}$, thus $\typ(S)=3$. By Proposition \ref{CM}, $\KK[S]$ is Cohen-Macaulay. Moreover, We have the equalities
\[
(5,10)+\a_1=(7,9)+\a_3=(8,8)+\a_5 \ , \ (8,8)+\a_2=(7,9)+\a_4=(5,10)+\a_3,
\]
which immediately imply that $S$ is nearly Gorenstein by Corollary \ref{NG}.
\end{example}

\begin{example} \label{ex2}
Let $\a_1=(2,0,0), \a_2=(0,2,0), \a_3=(0,0,2), \a_4=(1,1,0), \a_5=(1,0,1), \a_6=(0,1,1)$ and $S=\langle \a_1, \a_2, \a_3, \a_4, \a_5, \a_6 \rangle$. In this case $\Ap(S,E)=\{\0,\a_4,\a_5,\a_6\}$, and then $\typ(S)=3$. Using Proposition \ref{CM} it is possible to see that $\KK[S]$ is Cohen-Macaulay, whereas by the equalities 
$$\a_6+\a_1=\a_4+\a_5, \ \ \  \a_5+\a_2=\a_4+\a_6, \ \ \ \a_4+\a_3=\a_5+\a_6$$
and Corollary \ref{NG} it is easy to see that $S$ is nearly Gorenstein. 
\end{example}

\end{document}